\documentclass[final,narrowdisplay]{elsart1p} 
\usepackage{amssymb}
\usepackage{amsmath}
\usepackage[normal]{subfigure}
\include{TeXFigs} 
\usepackage{amsfonts}
\usepackage{amsmath} 
\usepackage{amssymb}
\usepackage{psfig}
\usepackage{verbatim}
\usepackage{amsmath}
\usepackage[mathcal]{euscript}
\usepackage{graphicx}
\numberwithin{equation}{section} 
\usepackage[mathcal]{euscript}
\newcommand \Xb {\overline X}
\newcommand \Ub {\overline U}

\newcommand \be         {\begin{equation}}
\newcommand \ee         {\end{equation}}

\newcommand \RR         {\mathbb{R}}
\newcommand \RN         {\mathbb{R}^N}

\newcommand \del        \partial
\newcommand \eps    \epsilon
\newcommand \lam    \lambda
\newcommand \gam    \gamma
\begin{document}

\begin{frontmatter}

\title{Why many theories of shock waves are necessary.
Convergence error in formally path-consistent schemes}

\author{Manuel J. Castro}

\address{Departamento de An\'alisis Matem\'atico,
Universidad de M\'alaga, 29071 M\'alaga, Spain.
\texttt{Castro@anamat.cie.uma.es}}

\author{Philippe G. LeFloch}

\address{Laboratoire J.-L. Lions \& Centre National de la Recherche Scientifique, 
\\
Universit\'e de Paris 6, 
4 Place Jussieu, 75252 Paris, France.
\texttt{LeFloch@ann.jussieu.fr}}

\author{Mar\'{\i}a Luz Mu\~noz-Ruiz}

\address{Departamento de Matem\'atica Aplicada,
Universidad de M\'alaga, 29071 M\'alaga, Spain. \texttt{Munoz@anamat.cie.uma.es}}

\author{Carlos Par\'es}

\address{Departamento de An\'alisis Matem\'atico,
Universidad de M\'alaga, 29071 M\'alaga, Spain.
\texttt{Pares@anamat.cie.uma.es}}

\date{}
%\today}

\begin{abstract} 
We are interested in nonlinear hyperbolic systems in nonconservative form arising in fluid dynamics,
and, for solutions containing shock waves, we investigate
the convergence of finite difference schemes applied to such systems. 
According to Dal Maso, LeFloch, and Murat's theory, 
a shock wave theory for a given nonconservative system requires prescribing a~priori 
a family of paths in the phase space. In the present paper, we consider schemes that 
are {\sl formally consistent}
with a given family of paths, and we investigate their limiting behavior as the mesh is refined. 
we first generalize to systems a property established earlier 
by Hou and LeFloch for scalar conservation laws, and we prove that nonconservative schemes 
generate, at the level of the limiting hyperbolic system, an {\sl convergence error source-term} which, 
provided the total variation of the approximations remains uniformly bounded, 
is a locally bounded measure. This convergence error measure is supported on the shock trajectories and, 
as we demonstrate here, is usually ``small''.
In the special case that the scheme converges in the sense of graphs
---a rather strong convergence property often violated in practice--- 
then this measure source-term vanishes. 
We also discuss the role of the equivalent equation associated with a difference scheme; here,
the distinction between scalar equations and systems appears most clearly since, for systems, 
the equivalent equation of a scheme that is formally path-consistent 
depends upon the prescribed family of paths. 
The core of this paper is devoted to 
investigate numerically the approximation of several (simplified or full) hyperbolic models 
arising in fluid dynamics. This leads us to the conclusion that 
for systems having 
nonconservative products associated with linearly degenerate characteristic fields,
the convergence error vanishes. For more general models, this measure 
is evaluated very accurately, especially 
by plotting the shock curves associated with each scheme under consideration; 
as we demonstrate, plotting the shock curves provide a convenient approach 
for evaluating the range of validity of a given scheme. 
\end{abstract}

\end{frontmatter}

 \tableofcontents

%=============================================================================

\section{Introduction}
\label{IN-0}

A number of non-conservative hyperbolic models have been introduced in fluid dynamics 
to serve as (simplified) models of two-phase or two-layer flows. Our objective in the present paper 
is to address the fundamental question whether finite difference schemes for nonconservative systems 
converge toward correct weak solutions containing shock waves. 
Addressing this important issue requires detailed numerical computations which we carry out here. 
The nonconservative hyperbolic systems under consideration have the general form
\be
u_t + A(u) \, u_x = 0, \qquad u=u(t,x) \in \RN,
\label{NONC}
\ee
where $u$ is the vector-unknown and $A=A(u)$ is a smooth, $N \times N$ matrix-valued map 
$A$ which admits real eigenvalues $\lam_1 < \ldots < \lam_N$ and a basis of eigenvectors $r_1, \ldots, r_N$. 
We are interested in solving the initial value problem associated with some initial condition
\be
u(0,x)=u_0(x), \qquad x \in \RR.
\label{NONC0}
\ee
The solutions of nonlinear hyperbolic systems are generally discontinuous; due to the non-divergence form of the equations the notion of solutions in the sense of distributions can not be used, and 
weak solutions satisfying \eqref{NONC} are defined in the sense introduced by 
LeFloch \cite{LeFloch0,LeFloch1,PLF-1,LeFloch2,PLF-graph} 
and Dal~Maso, LeFloch, and Murat \cite{DLM} (cf.~Section~2 below for 
a brief review of the theory).

Generally speaking, solutions to \eqref{NONC} depend upon regularization mechanisms;
for instance, different approximation schemes may converge toward different solutions, and for this reason 
in developing the well-posedness theory, higher-order regularization effects such as 
viscosity, capillarity, relaxation terms, must be taken into account in the modeling. 
For instance, for continuous models one may consider the regularization
\be
u_t^\eps + A(u^\eps) \, u_x^\eps = R^\eps,
\label{NONC2}
\ee
where $R^\eps$ depends upon higher-derivatives of $u$ together with (one or several) small-scale parameter(s)
$\eps$; the physically meaningful solutions are defined as the singular limits 
$$
u := \lim_{\eps \to 0} u^\eps. 
$$
Furthermore, as established in \cite{LeFloch1}, shock waves in such solutions are determined 
by traveling wave solutions to \eqref{NONC2}, that is, Rankine-Hugoniot relations for shock waves 
are determined from the given regularization. 

In the present paper, we demonstrate that while, for certain simplified models, 
solutions are actually stable upon regularization
and the detailed knowledge of the right-hand side of \eqref{NONC2} is unnecessary,
however for general systems such as the ``full'' systems of two-phase flows,
the general DLM theory is necessary. 
Still, as pointed out by Hou and LeFloch \cite{HouLeFloch} ---who focused
attention on the same issues for nonconservative formulations of scalar hyperbolic equations---
and by Hayes and LeFloch \cite{HayesLeFloch0,HayesLeFloch} and LeFloch and Mohammadian \cite{LM} 
---who studied the effect of diffusive and dispersive terms---
the effects of the regularization $R^\eps$ may be difficult to pinpoint in practice. 
In view of the fact that the models under study are derived from modeling approximation assumptions, 
this fully justifies the use of a numerical strategy based on a direct discretization of the nonconservative 
hyperbolic models \eqref{NONC}. Our conclusions, therefore, justify to search for
robust and efficient high-order schemes for the approximation of nonconservative systems.
In particular, Berthon and Coquel \cite{BC,BerthonCoquel} and Chalons and Coquel \cite{CC}, 
have introduced various numerical strategies for models of complex fluid flows including turbulence models, 
while Par\'es \cite{Pares} and Mu\~noz-Ruiz and Par\'es \cite{MP} have developed many important applications. 

Finally, for anther standpoint to the theory and the numerical analysis of nonconservative products, 
we refer to Berthon, Coquel, and LeFloch \cite{BCL} who connected the theory of nonconservative products 
with the concept of a kinetic relation \cite{LeFloch-book}. They 
introduced a general framework to handle nonconservative systems; 
this framework encompasses a large number of examples arising in the 
applications. In particular, they rigorously analyzed a typical model of turbulent fluid dynamics 
by establishing the existence and properties of a physically relevant family of traveling waves 
and deriving the corresponding kinetic function.

%*****************************************************************************************

\section{The convergence error measure}
\label{Conver}

\subsection{DLM familes of paths and nonconservative products}

Let $\Omega$ be an open subset of $\RN$ and $g : \Omega \to \Omega$ be a smooth mapping.
Given a function with bounded variation $u: \RR \to \Omega$,
the definition introduced by Dal Maso, LeFloch, and Murat \cite{DLM}
allows one to define products of the form $g(u) \, {du \over dx}$
provided a family of Lipschitz continuous paths $\Phi: [0,1] \times \Omega \times \Omega$ is prescribed,
which must satisfy certain natural regularity conditions, in particular
\be \label{cond1}
\Phi(0;u_l, u_r) = u_l, \qquad  \Phi(1;u_l, u_r) = u_r,
\ee
\be \label{cond2}
\Phi(s;u_l, u_l) = u_l,
\ee
for all $s \in [0,1]$ and $u_l, u_r \in \Omega$. The nonconservative product, denoted
by $\big[g(u) \, {du \over dx}\big]_\Phi$, is defined as a bounded measure, which is absolutely continuous with respect
to the total variation measure of the function $u$ and, in particular,
coincides with the distributional derivative ${d \over dx} f$ in the special case of a conservative product, for which $g(u) = Df(u)$ for some $f$.

The DLM theory was applied to nonconservative systems of the form \eqref{NONC};
the Riemann problem was solved and, later, the general Cauchy problem \cite{PLF-Liu}. 
In the course of this analysis, the notion of $\Phi$-completion
$(X,U_\Phi)$ of the graph of a BV function $u$ was introduced. The key stability result
in \cite{DLM} is the following: if $u^\Delta:\RR \to \Omega$ is a sequence of BV functions with uniformly
bounded amplitude and total variation
\be
\label{bd}
\sup_\RR |u^\Delta| + TV_\RR(u^\Delta) \lesssim 1,
\ee
converging almost everywhere to a limit function $u:\RR \to \Omega$, then a sufficient condition for the corresponding nonconservative products to converge
\be
\label{conve}
\Big[g(u^\Delta) \, {du^\Delta \over dx}\Big]_\Phi \rightharpoonup \Big[g(u) \, {du \over dx}\Big]_\Phi
\ee
in the weak-star sense of measures,
is that their $\Phi$-completions  $(X^\Delta,U_\Phi^\Delta)$ converge in the uniform distance of graphs
precisely to the $\Phi$-completion $(X,U_\Phi)$ of the limit $u$.

%--------------------------------------------------------------------------------------------------------------------

\subsection{A class of finite difference schemes}
\label{schemes}

For the approximation of nonconservative systems we introduce here a general family of numerical schemes which includes, in particular, three classes of schemes
of particular interest: Godunov, Roe, and Lax-Friedrichs. 

For the discretization of the initial value problem \eqref{NONC}-\eqref{NONC0}, 
we introduce computing cells $I_i=[x_{i-1/2},x_{i+1/2}]$ and, for simplicity in the presentation, we
assume that these cells have constant size $\Delta=\Delta x$. 
We also define $x_{i+1/2}=i\Delta x$ and $x_i=(i-1/2)\Delta x$, the latter being the center
of the cell $I_i$. Finally we denote by $\Delta t$ the (constant) time length and we set $t^n=n\Delta t$.
We denote by $u_i^n$ the approximation of the cell averages of the exact solution provided by the numerical scheme:
$$
u_i^n\cong\frac{1}{\Delta x}\int_{x_{i-1/2}}^{x_{i+1/2}}u(t^n,x)\,dx.
$$

We are interested in schemes of the general form
\be
\label{scheme}
u_i^{n+1} = u_i^n - \frac{\Delta t}{\Delta x} \bigl( M_{i-1/2}^{n,+} + M_{i+1/2}^{n,-}\bigr),
\ee
where
$$
M_{i+1/2}^{n,\pm }= M^\pm(u^n_{i-q}, \dots, u^n_{i+p}).
$$

From now on a DLM family of paths $\Phi$ for the nonconservative system \eqref{NONC} is fixed. 
Following Par\'es \cite{Pares}, we consider {\sl formally path-consistent schemes,} i.e. schemes that 
are consistent with the family of paths $\Phi$ in the following sense:
$M^-$ and $M^+$ are Lipschitz continuous mappings from $\Omega^{p+q+1}$ to $\Omega$ satisfying:
\be
\label{consist}
M^\pm(u, \dots, u) = 0, \qquad \ u \in \Omega,
\ee
and for every $u_i \in \Omega$, $i = -q, \dots, p$,
\be 
\label{jumps}
M^-(u_{-q}, \dots, u_{p} ) + M^+ (u_{-q}, \dots, u_{p} ) = \int_0^1A(\Phi(s;u_0,u_1))\frac{\partial\Phi}{\partial s}(s;u_0,u_1)\, ds.
\ee
These conditions provide a generalization to
the concept of conservative scheme introduced by Lax for systems of conservation laws and, for this 
reason, were originally called ``path-conservative'' in \cite{Pares}. 
As we will see, this definition --although quite natural-- needs to be handled carefully.

It is convenient to assume some particular structure on the given family of paths. Precisely, we assume that the matrix 
$A$ and the paths $\Phi$ satisfy the following restrictions: 
\begin{itemize}

\item[(R1)] Given an integral curve $\gamma$ of a linearly
degenerate field and $u_l, u_r \in \gamma$, the path
$\Phi(s;u_l,u_r)$ is a parametrization of the arc of $\gamma$ connecting ${u}_l$ and ${u}_r$.

\

\item[(R2)] Given an integral curve $\gamma$ of a  genuinely
nonlinear field  and $u_l, u_r \in \gamma$, with  $\lambda(u_l) < \lambda(u_r)$,  being
$\lambda(u)$ the corresponding eigenvalue, the path $\Phi(s;u_l,u_r)$ is a
parame\-trization of the arc of $\gamma$ connecting ${u}_l$ and
${u}_r$.

\

\item[(R3)] Let us denote by $ {\mathcal{RP}} \subset \Omega
\times \Omega$ the set of pairs $({u}_l, {u}_r)$ for which the
Riemann problem
\begin{equation}\label{RP}
\begin{cases}
\displaystyle{ u_t + A(u)u_x = 0,}   \smallskip \\
{u}(x,0)=\begin{cases}
        {u}_l, & x<0, \\
        {u}_r, &  x>0,
        \end{cases}
\end{cases}
\end{equation}
has a unique self-similar weak solution composed by $N$ (possibly trivial) 
simple waves connecting $N+1$ intermediate constant states
$$
u_0 = u_l, \, u_1, \dots, u_{N-1}, \, u_N = u_r.
$$
Given $({u}_l, {u}_r) \in  {\mathcal{RP}}$,
the curve described by the path ${\Phi }(\cdot; {u}_l, {u}_r)$ is
equal to the union of those corresponding to the paths ${\Phi
}(\cdot; {u}_{j-1}, {u}_j)$, $j=1, \dots, N$.

\end{itemize}

\

We now present several schemes of particular interest in the present paper. 
First of all, the Godunov method for the nonconservative system \eqref{NONC} takes 
the form  \eqref{scheme} with
$$
\aligned 
& M_{i+1/2}^{n,-} & = & \int_0^1A(\Phi(s;u_i^n,u_{i+1/2}^n))\frac{\partial\Phi}{\partial s}(s;u_i^n,u_{i+1/2}^n)\, ds, 
\\
& M_{i+1/2}^{n,+} & = & \int_0^1A(\Phi(s,u_{i+1/2}^n,u_{i+1}^n))\frac{\partial\Phi}{\partial s}(s;u_{i+1/2}^n,u_{i+1}^n)\,ds, 
\endaligned
$$
where $u_{i+1/2}^n$ is the (constant) value at $x =0$ of the solution of the Riemann problem consisting of  \eqref{NONC} with
initial condition:
$$
u(x,0) = \begin{cases}
u_i^n,       &  x < 0,
\\
u_{i+1}^n,   & x >0.
\end{cases}
$$
In the derivation it is important to assume a ``CFL-1/2 condition'', as noted in \cite{MP}.  

In the case in which the solution of the Riemann problem is discontinuous at $x=0$, such discontinuity has to be stationary and we can replace  $u_{i+1/2}^n$ either by  the limit of the solution to the left or
to the right of $0$.

Secondly, Roe methods provide linear approximate Riemann solvers. In view of the framework in \cite{DLM}
the following generalization to Roe's standard approach was proposed by LeFloch \cite{LeFloch2}. 
Given a family of paths $\Phi$, a function $A_\Phi\colon\Omega\times\Omega\mapsto\mathcal{M}_{N\times
N}(\mathbb{R})$ is called a (generalized) Roe linearization if for all $u_l, u_r \in \Omega$ the 
following properties hold: 
\begin{enumerate}
\item $A_\Phi(u_l,u_r)$ has $N$ distinct real eigenvalues,
\item $A_\Phi(u_l,u_l)=A(u_l)$,
\item $A_\Phi(u_l,u_r)\cdot (u_r-u_l)= \displaystyle \int_0^1A(\Phi(s;u_l,u_r))\frac{\partial\Phi}{\partial
s}(s;u_l,u_r)\,ds$.
\end{enumerate}

Once a Roe linearization has been chosen, the corresponding Roe
scheme takes the form \eqref{scheme} with
$$
\aligned 
& M_{i+1/2}^{n,-} = A^{n,-}_{i+1/2}\cdot (u_{i+1}^n-u_i^n), 
\\
& M_{i+1/2}^{n,+} = A^{n,+}_{i+1/2}\cdot (u_{i+1}^n-u_i^n), 
\endaligned 
$$
where
$A^n_{i+1/2}=A_\Phi(u_i^n,u_{i+1}^n)$, 
$$
{\mathcal{L}}_{i+1/2}^{n,\pm }=\left[
\begin{array}{ccc}
(\lambda_{i+1/2,1}^n)^{\pm } &  & 0 \\
& \ddots  &  \\
0 &  & (\lambda_{i+1/2,N}^n)^{\pm }
\end{array}
\right],
$$
and 
$$
{A}_{i+1/2}^{n,\pm }=
      \mathcal{K}^n_{i+1/2}{\mathcal{L}}_{i+1/2}^{n,\pm}\left(\mathcal{K}^n_{i+1/2}\right)^{-1}. 
$$
Here, ${\mathcal{L}}^n_{i+1/2}$ denotes the diagonal matrix whose
coefficients are the eigenvalues of $ A^n_{i+1/2}$, $\lambda^n_{i+1/2,1}<\lambda^n_{i+1/2,2}<\cdots
<\lambda^n_{i+1/2,N}$, and $\mathcal{K}^n_{i+1/2} $ a  $N\times N$ matrix  whose columns
are associated eigenvectors.

\

Third, a generalization of the classical Lax-Friedrichs method  to \eqref{NONC}  is given by \eqref{scheme} with
the choice:
\be
\label{DLFM} 
\aligned 
M_{i+1/2}^{n,-} & = \int_0^1\widehat{A}^-(\Phi(s;u_i^n,u_{i+1}^n))\frac{\partial\Phi}{\partial s}(s;u_i^n,u_{i+1}^n)\, ds, \\
M_{i+1/2}^{n,+} & = \int_0^1\widehat{A}^+(\Phi(s,u_{i}^n,u_{i+1}^n))\frac{\partial\Phi}{\partial s}(s;u_{i}^n,u_{i+1}^n)\,ds,
\endaligned
\ee 
where
$$
\widehat{A}^{-}(u) = \frac{1}{2}\left(- \frac{\Delta x}{\Delta t} Id + A(u)\right)
\label{Alf-}, 
\qquad 
\widehat{A}^{+}(u) = \frac{1}{2}\left(\frac{\Delta x}{\Delta t}
Id +A(u)\right), 
$$
being $Id$ the $N\times N$ identity matrix.  

%-----------------------------------------------------------------------------------------------------

\subsection{Convergence to a nonconservative system with measure-source term}

We denote by $u^\Delta$ the sequence of piecewise constant approximate solutions generated by 
a finite difference scheme of the form described above. 
By extending the arguments in Hou and LeFloch \cite{HouLeFloch} we can prove:

\

\begin{claim}
\label{conv}
Consider a nonconservative hyperbolic system \eqref{NONC} together
with a given family of paths $\Phi$. Suppose that $u^\Delta$ is a sequence of approximate solutions constructed by one of
the finite difference schemes described in Subsection~\ref{schemes}
and satisfying the bounds \eqref{bd} uniformly in time. Suppose that the scheme is
{\rm formally consistent with the family of paths} 
$\Phi$. 
Then, given any subsequence of $u^\Delta$ converging almost everywhere to some limit,
denoted by $v$, the following holds: 
\begin{enumerate}
\item There exists a bounded measure $\mu_v : \RR_+ \times \RR \to \RN$ 
(called the {\rm convergence error measure})
such that the limit $v$ satisfies the following hyperbolic system with source-term
\be
\label{source}
v_t + \big[A(v) \, v_x\big]_\Phi = \mu_v.
\ee
\item 
Moreover, when the $\Phi$-completion of the graphs of $u^\Delta$ converges in the uniform 
sense of graphs
towards the $\Phi$-completion of the limit $v$, i.e. 
$$
(X^\Delta, U^\Delta_\Phi) \to (Y,V_\Phi), 
$$ 
then the convergence error measure $\mu_v$ vanishes
identically and $v$ is a weak solution to the system
\be
\label{exact}
v_t + \big[A(v) \, v_x\big]_\Phi = 0.
\ee
\end{enumerate} 
\end{claim}

\

The above result can be interpreted as a ``nonconservative extension'' to the classical Lax-Wendroff 
theorem for systems of conservation laws \cite{LaxWendroff}. 
It should be observed that the convergence in the sense of graphs is very strong; it does hold
for the Glimm and front tracking schemes, but usually fails for finite difference schemes.

Our main objective in the present paper is to investigate the source and the amplitude of this convergence error,  
which can be measured in terms of the measure $\mu_v$ or, equivalently,  
in  terms of the Rankine-Hugoniot curves associated with the given scheme. Indeed, computing numerically 
the shock curves associated with various schemes of interest is one of the main purposes of this work.

\ 
 
\noindent{\sl Proof.} We follow Hou and LeFloch \cite{HouLeFloch} and decompose the scheme into 
a part that converges to the hyperbolic system and an error term, 
and we then rely on stability results established in the general DLM theory \cite{DLM}. 
We need to prove that
\be
\label{weak_form}
\int_0^{\infty} \int_{-\infty }^{\infty } v(t,x)\varphi_t(t,x)\,dxdt
-\int_0^{\infty } \big\langle \big[ A(v(t, \cdot ))v_x(t, \cdot ) \big]_{\Phi }, \varphi (t, \cdot ) \big\rangle \, dt = 0
\ee
for all compactly supported test-function $\varphi = \varphi(t,x)$.

We set 
$$
\varphi_i^n=\varphi (t^n,x_i), \qquad \varphi_{i+1/2}^n = \varphi (t^n,x_{i+1/2}) 
$$
and we multiply the discrete equation \eqref{scheme} by $\varphi_i^n$.
After summing over $i$ and $n$, and then applying summation by parts we obtain the identity 
\be
\label{discrete_weak_form}
\aligned
& \Delta x\Delta t \sum_{n=1}^{\infty } \sum_{i=-\infty }^{\infty } u_i^n \frac{\varphi_i^n-\varphi_i^{n-1}}{\Delta t} \\
& - \Delta t  \sum_{n=0}^{\infty } \sum_{i=-\infty }^{\infty }  \big( M_{i+1/2}^{n,-} + M_{i+1/2}^{n,+} \big) \big( \varphi_{i+1/2}^n + O(\Delta x) \big) =0.
\endaligned
\ee

We want to prove that equation \eqref{weak_form} can be obtained by passing to the limit as $\Delta$ tends to 0 in \eqref{discrete_weak_form}. The convergence of the first term in \eqref{discrete_weak_form} to the corresponding one in \eqref{weak_form}
is obtained as in the conservative case. Concerning the second term, since $M^{\pm}$ are Lipschitz 
continuous and the states $u_i^n$ are uniformly bounded, we only have to study the convergence of
\begin{multline*}
\Delta t \sum_{n=0}^{\infty } \sum_{i=-\infty }^{\infty }  \big( M_{i+1/2}^{n,-} + M_{i+1/2}^{n,+} \big) \varphi_{i+1/2}^n \\
=\Delta t \sum_{n=0}^{\infty } \sum_{i=-\infty }^{\infty }  \left( \int_0^1 A  \left( \Phi (s;u_i^n,u_{i+1}^n) \right) \Phi_s (s;u_i^n,u_{i+1}^n) \,ds \right) \varphi_{i+1/2}^n,
\end{multline*}
where we have used our assumption that the scheme is formally consistent with the given family of paths $\Phi$. 

Using the definition of the product $A(u^\Delta(t, \cdot ))u^\Delta_x(t, \cdot )$ as a measure this term can be rewritten as follows:
$$
\aligned
& \int_0^{\infty } \big\langle \big[ A(u^\Delta(t, \cdot )) u^\Delta_x (t, \cdot ) \big]_{\Phi }, \varphi (t, \cdot ) \big\rangle \, dt 
\\
& + \sum_{n=0}^{\infty } \sum_{i=-\infty }^{\infty } \int_{t_n}^{t_{n+1}} \left( \int_0^1 A  \left( \Phi (s;u_i^n,u_{i+1}^n) \right) \Phi_s (s;u_i^n,u_{i+1}^n) \, ds \right) \big( \varphi_{i+1/2}^n -\varphi (t, x_{i+1/2}) \big) \, dt.
\endaligned
$$
The second summand trivially converges to $0$. The proof is concluded by using the main stability result 
and error estimate in~\cite{DLM}. Indeed, if the sequence $u^\Delta$ converges in the sense of graphs, 
then, as recalled in \eqref{conve} the nonconservative product converges in the weak-star sense of measures,
$$
\int_0^{\infty } \big\langle \big[ A(u^\Delta(t, \cdot )) u^\Delta_x (t, \cdot ) \big]_{\Phi }, \varphi (t, \cdot ) \big\rangle \,dt 
\longrightarrow
\int_0^{\infty } \big\langle \big[ A(v(t, \cdot )) v_x (t, \cdot ) \big]_{\Phi }, \varphi (t, \cdot ) \big\rangle \,dt 
$$
and we recover the exact system \eqref{exact}. In the general case, the $\Phi$-graph completion 
$(X^\Delta, U^\Delta)$ of 
$u^\Delta$ converges to some limiting graph $(\Xb, \Ub)$, 
$$
(X^\Delta, U^\Delta) \to (\Xb, \Ub), 
$$ 
which is such that its projected BV function coincides with the pointwise limit $v$ of $u^\Delta$. 
In turn the limit nonconservative product is based on the BV function $v$ and we obtain 
$$
\aligned
& \int_0^{\infty } \big\langle \big[ A(u^\Delta(t, \cdot )) u^\Delta_x (t, \cdot ) \big]_{\Phi }, \varphi (t, \cdot ) \big\rangle \,dt 
\\
& \longrightarrow
\int_0^\infty \int_\RR A(\Ub(t, \cdot )) \Ub (t, \cdot ) \varphi (t, \cdot ) \big\rangle \,dt 
=
\int_0^\infty  
\big\langle \big[ A(v(t, \cdot)) v_x (t,\cdot) \big]_\Phi + \mu_v(t), \varphi (t, \cdot ) \big\rangle \, dt 
\endaligned
$$
for some time-dependent measure $\mu_v$. This completes the proof. \qed

%*****************************************************************************************

\section{Equivalent equations for nonconservative systems}
\label{Equiv}

\subsection{Derivation of the equivalent equation}

In this section, we  derive the equivalent equations corresponding to the Lax-Friedrichs-type scheme introduced in Subsection~\ref{schemes}, which reads 
$$
\aligned
\frac{1}{\Delta t}\left( u_i^{n+1}-\frac{1}{2}\left( u_{i-1}^n+u_{i+1}^n \right) \right) 
+ \frac{1}{2\Delta x}\Big( 
& \int_0^1 A \left( \Phi (s;u_{i-1}^n,u_i^n) \right) \frac{\partial \Phi }{\partial s} (s;u_{i-1}^n,u_i^n) \, ds  
\\
& + \int_0^1 A \left( \Phi (s;u_i^n,u_{i+1}^n) \right) \frac{\partial \Phi }{\partial s} (s;u_i^n,u_{i+1}^n) \, ds \Big) = 0. 
\endaligned
$$
We consider first the equivalent equation at second order. Performing a formal Taylor expansion in the scheme we obtain:
$$
\aligned
& v_t + A(v)v_x + \frac{\Delta t}{2} \left( v_{tt} - \frac{\Delta x^2}{\Delta t^2}\, v_{xx} + \frac{\Delta x}{\Delta t} \, I_2(v) \right) = 0,
\\
& I_2(v) = \int_0^1 DA\left( v \right) \left( D_{u_l} \Phi \cdot v_x, D_{u_l} \Phi_s \cdot v_x \right) \, ds 
+ \int_0^1 DA\left( v \right) \left( D_{u_r} \Phi \cdot v_x, D_{u_r} \Phi_s \cdot v_x \right) \, ds,
\endaligned
$$
where the following notation has been used:  given two vectors
$v = [v_1, \dots, v_N]^T$, $w= [w_1, \dots, w_N]^T$, 
$D{A}(u)(v,w)$ represents the derivative of ${A}(u)$ in the direction of the vector $v$ (which is a
matrix) applied to the vector $w$, i.e.
$D{A}(u)(v,w) = \left(\sum_{k=1}^N v_k
\partial_{u_k}{A}(u) \right) \cdot w$, 
where $\partial_{u_k}{A}(u)$ is the $N \times N$ matrix
whose $(i,j)$ entry is $\partial_{u_k}a_{i,j}(u)$.

As usual, we write the $t$-derivatives using $x$-derivatives in order to obtain a new modified equation.
Using the equality
$$
\left( A^2(v)v_x \right)_x = DA(v) \left( v_x, A(v)v_x \right) + A(v) \left( A(v)v_x \right)_x,
$$
we obtain the modified equation for the Lax-Friedrichs method:
\begin{multline} \label{modified_2}
v_t + A(v)v_x = \frac{\Delta x^2}{ 2\Delta t} \left( v_{xx} - \frac{\Delta t^2}{\Delta x^2} \left( A^2(v)  v_x \right)_x \right. \\
\left. -  \frac{\Delta t^2}{\Delta x^2}  \left( DA(v) \left( A(v)v_x, v_x \right) - DA(v) \left( v_x, A(v)v_x \right) \right) \right)- \frac{\Delta x}{2}I_2(v) .
\end{multline}
Note that, in the expression of the modified equations, the only term that depends on the choice of the family of paths
is the last one: $\Phi$ only appears in the expression of $I_2(v)$.

In a similar way, by a simple (but tedious) calculation we can obtain the equivalent equations at third order: 
$$
\aligned 
& v_t + A(v)v_x 
\\
= \, 
& \Delta x \left( \frac{\Delta x}{2\Delta t}v_{xx} - \frac{\Delta t}{2\Delta x} \Theta (v) - \frac{1}{2}I_2(v) \right) 
   + \Delta x^2 \left( \frac{1}{2} \left( A(v)v_{xx} \right)_x + \frac{1}{6} A(v)v_{xxx} \right. 
\\
&- \frac{\Delta t^2}{3\Delta x^2} \left( \left( A(v)\Theta (v) \right)_x - DA(v)\left( \Theta (v), v_x\right) - DA(v)\left( v_x,\Theta (v)\right)  \right) 
\\
&- \frac{\Delta t^2}{12\Delta x^2} D^2A(v)\left( A(v)v_x, A(v)v_x, v_x \right) - \frac{\Delta t^2}{6\Delta x^2} DA(v)\left( A(v)v_x, \left(A(v)v_x\right)_x \right) 
\\
&- \frac{\Delta t}{4\Delta x} \left( \left( A(v)I_2(v) \right)_x + DA(v) (I_2(v), v_x) + DA(v) (v_x,I_2(v)) \right) 
\\
&\left. - \frac{1}{4} D^2A(v)(v_x,v_x,v_x) - \frac{1}{4}I_3(v) - \frac{\Delta t}{4\Delta x} \left( I_2(v) \right)_t \right),
\endaligned
$$
where we have used the notation
$$
D^2{A}(u)(v^1,v^2,w) = \left(\sum_{k,m=1}^N v_k^1 v_m^2
\partial_{u_ku_m}^2{A}(u) \right) \cdot w
$$
and
$$
\Theta (v) = \left( A^2(v)v_x \right)_x + DA(v) \left( A(v)v_x, v_x \right) - DA(v) \left( v_x, A(v)v_x \right). 
$$

Furthermore, the expansion of $\left( I_2(v) \right)_t$ is
$$
\left( I_2(v) \right)_t =  I_{2,1}(v) + I_{2,2}(v) + I_{2,3}(v), 
$$
with 
$$
\aligned
I_{2,1}(v) 
=  & \int_0^1 D^2A(v) \left( A(v)v_x, D_{u_l} \Phi \cdot v_x, D_{u_l} \Phi_s \cdot v_x \right) \, ds 
\\
& + \int_0^1 D^2A(v) \left( A(v)v_x, D_{u_r} \Phi \cdot v_x, D_{u_r} \Phi_s \cdot v_x \right) \, ds 
\\
& + \int_0^1 DA(v) \left( D^2_{u_l u_l} \Phi \left( A(v)v_x,v_x \right), D_{u_l} \Phi_s \cdot v_x \right) \, ds 
\\
& + \int_0^1 DA(v) \left( D^2_{v_r v_r} \Phi \left( A(v)v_x,v_x \right), D_{u_r} \Phi_s \cdot v_x \right) \, ds, 
\endaligned
$$

$$
\aligned
I_{2,2}(v) 
=  
& \int_0^1 DA(v) \left( D_{u_l} \Phi \cdot \left( DA(v) (v_x,v_x) + A(v)v_{xx} \right), D_{u_l} \Phi_s \cdot v_x  \right) \, ds 
\\
& + \int_0^1 DA(v) \left( D_{u_r} \Phi \cdot \left( DA(v) (v_x,v_x) + A(v)v_{xx} \right), D_{u_r} \Phi_s \cdot v_x  \right) \, ds 
\\
& + \int_0^1 DA(v) \left( D_{u_l} \Phi \cdot v_x, D^2_{u_l u_l} \Phi_s \left( A(v)v_x,v_x \right) \right) \, ds 
\\
& + \int_0^1 DA(v) \left( D_{u_r} \Phi \cdot v_x, D^2_{u_r u_r} \Phi_s \left( A(v)v_x,v_x \right) \right) \, ds, 
\endaligned
$$
and 
$$
\aligned
I_{2,3}(v) 
=  & \int_0^1 DA(v) \left( D_{u_l} \Phi \cdot v_x, D_{u_l} \Phi_s \cdot \left( DA(v) (v_x,v_x) + A(v)v_{xx} \right) \right) \, ds 
\\
& + \int_0^1 DA(v) \left( D_{u_r} \Phi \cdot v_x, D_{u_r} \Phi_s \cdot \left( DA(v) (v_x,v_x) + A(v)v_{xx} \right) \right) \, ds.
\endaligned
$$
The above formula illustrate the high complexity of the equivalent equation approach when dealing
with nonconservative schemes. 

%----------------------------------------------------------------------------------------------------------------

\subsection{Role of the equivalent equation}

These modified equations are useful to understand why the numerical solutions may not converge to the weak
solutions of the system. Let us suppose that the family of paths $\Phi$ used in the definition of weak solutions
is based on a parabolic regularization of the system
\be
u_t^\eps + A(u^\eps) \, u_x^\eps = \eps (D(u^\eps)u^\eps_x)_x,
\label{NONC-p}
\ee
where $D$ is a viscosity matrix which is admissible in the following sense:  if $u_l$, $u_r$ can be connected by
a discontinuity satisfying the Rankine-Hugoniot conditions related to the family of paths:
\be
\xi (u_r - u_l) = \int_0^1 A(\Phi(s; u_l, u_r)) \frac{\partial \Phi}{\partial s}(s; u_l, u_r) \, dx
\label{RH}
\ee
for some $\xi \in \RR$, then $\Phi(s; u_l, u_r)$ is a reparametrization of the solution of the differential
system:
\be
-\xi v' + A(v)v' = (D(v)v')',
\ee
with the conditions:
\be
\lim_{\xi \to -\infty}v(\xi) = u_l, \quad \lim_{\xi \to \infty}v(\xi) = u_r.
\ee
Under these hypotheses it can be checked \cite{LeFloch1} that the function
\be
u(x,t) = \begin{cases}
u_l & x < \xi t,\\
u_r & x > \xi t.
\end{cases}
\ee
is a weak solution in the sense of Dal~Maso, LeFloch, and Murat. 

If now the Lax-Friedrichs scheme is applied to the hyperbolic system, the limits of the numerical solutions
provide approximations to the vanishing viscosity limits related to the regularization
\eqref{modified_2}  which is different of \eqref{NONC-p}. The difficulty comes from
the fact that, unlike the conservative case, the vanishing viscosity limits depend on the regularization
of the problem.
Even if, for simplicity, we have only calculated the modified equations corresponding to the Lax-Friedrichs
scheme, the same difficulty would be present for any other scheme involving a numerical viscous
term: the numerical solutions approximate the vanishing viscosity limit of a modified equation whose
regularization terms depend both on the chosen family of paths and on the specific form of its viscous terms.

%*****************************************************************************************

\section{Examples of nonconservative hyperbolic systems}

\subsection{A simplified model}
\label{csimple}

We begin with a hyperbolic system containing nonconservative products,
which will be used in the following section to perform numerical experiments.
We consider the system  
\begin{equation}
\aligned 
& h_t + q_x = 0,
\\
& q_t + \left(\frac{q^2}{h}\right)_x + qhh_x = 0,
\endaligned
\label{simple}
\end{equation}
which has the form $w_t + A(w)w_x = 0$ with
$$
w = \left[ \begin{array}{c} h \\ q \end{array} \right],\quad
 A(w)
= \left[
\begin{array}{cc}
0 & 1 \\
-u^2 + u^2h & 2u
\end{array}
\right],
$$
and $u = q/h$. In the region
$$ 
\Omega = \{ (h,q) \, | \, 0<q, \, 0 < h < (16 q)^{1/3} \},  
$$
the system is strictly hyperbolic and all characteristic
fields are genuinely nonlinear. The eigenvalues of this system are
$$
\lambda_1 = u - h\sqrt{u}, \quad \lambda_2 = u + h\sqrt{u},
$$
and the integral curves of the first  and second characteristic fields
are 
$$ 
\sqrt{u} + h/2 = const, \qquad \sqrt{u} - h/2 = const,
$$
respectively.

In order to define the jump conditions, the paths connecting the left-
and right-hand limits $w^\pm = [h^\pm, q^\pm]$ at a shock are chosen to be
the union of the segment connecting $w^-$ with $w^* = [h^+, q^-]$ and
the segment connecting $w^*$ to $w^+$:
\begin{equation}\label{camino_simple}
\Phi (s; w^-, w^+)
= \left\{
\begin{array}{cc}
\left[ \begin{array}{c} h^- + 2s(h^+ -h^-) \\ q^-
\end{array} \right], & \qquad \mbox{$0 \leq s \leq 1/2$;}\\
  & \\
\left[ \begin{array}{c} h^+ \\ q^- + (2s -1 )(q^+ - q^-)
\end{array} \right],
  & \qquad \mbox{$1/2 \leq s \leq 1$;}  \end{array} \right.
  \end{equation}
Such paths were introduced in \cite{DLM} to illustrate certain issues encountered with nonconservative products. 
The  corresponding jump conditions \eqref{RH} are the following:
$$
\aligned 
& \xi[h] = [q],
\\
& \xi[q] = \left[ \frac{q^2}{h}\right] + q^-\left[
\frac{h^2}{2} \right].
\endaligned
$$

\begin{center}
\begin{figure}[h!ptb]
\centering
\includegraphics[width=10.5cm,height=10.0cm]{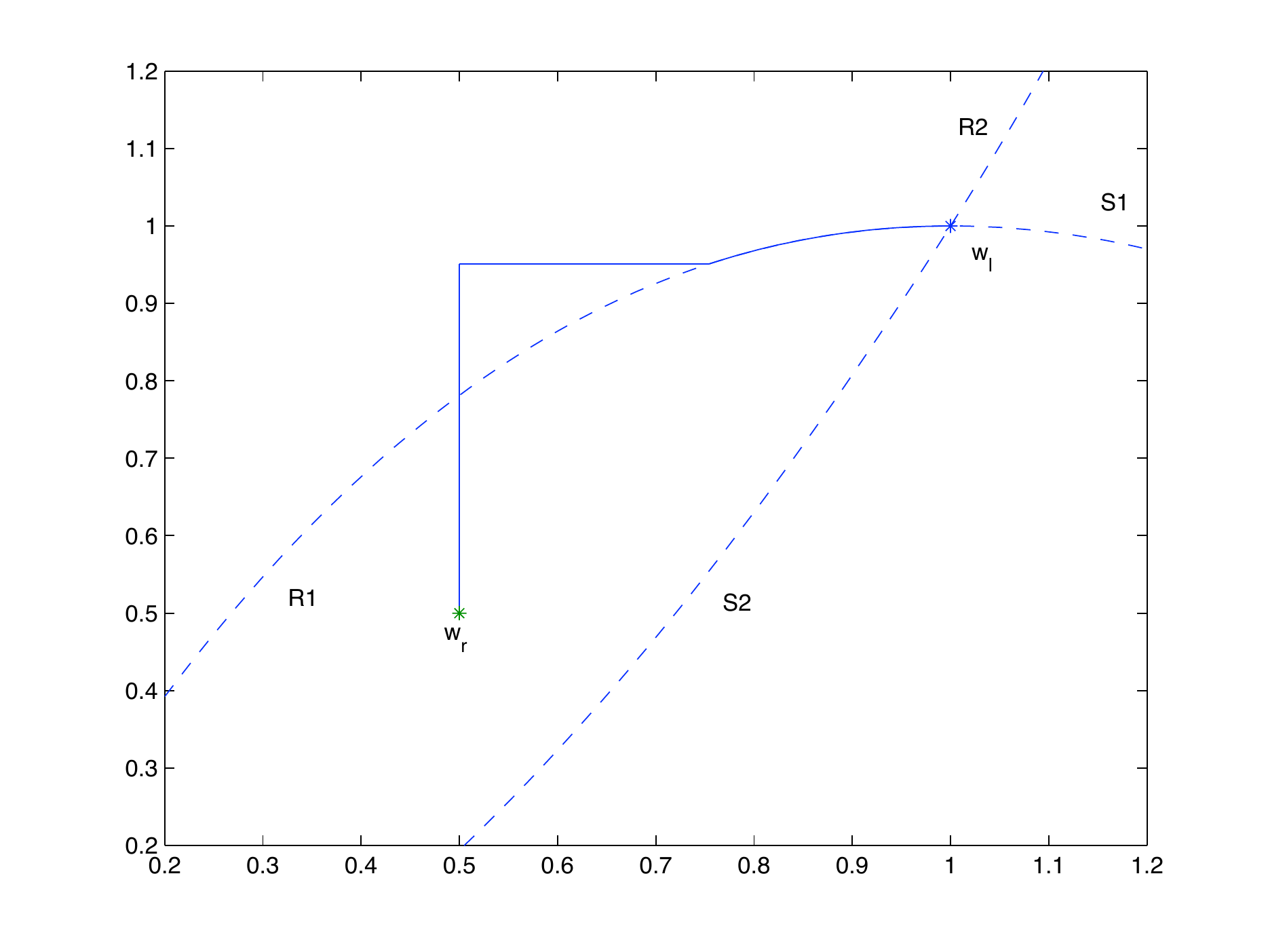}
\caption{Path connecting the states $w_l$ and $w_r$ (continuous line) and
set of states that can be connected to $w_l$ by an entropy satisfying 1 or 2-wave
(dashed lines).}
\label{figpath}
\end{figure}
\end{center}

Once the jump conditions have been stated, it is possible to solve the
Riemann problem for any pair of states which are sufficiently enough. 
Then, a family of paths
satisfying (R2) and (R3) is constructed.  In Figure~\ref{figpath} the path linking the states $w_l = [1,1]^T$ and
$w_r = [0.5, 0.5]^T$ is depicted.  In this case, the  solution
of the Riemann problem consists of a $1$-rarefaction connecting $w_l$ to
an intermediate state $\tilde{w}$ and a $2$-shock linking $\tilde{w}$ to $w_r$.
As a consequence, the path consists of the arc of the 1-integral curve linking $w_l$ to
$\tilde{w}$ and two segments connecting $\tilde{w}$ and $w_r$ which are parallel to the axis .
The set of states that can be connected to $w_l$  by an entropy satisfying
1-wave or 2-wave are also depicted. In each case $Ri$ denotes $i$-rarefactions
and $Si$ denotes $i$-shocks. 

%-------------------------------------------------------------------------------------------------------

\subsection{A class of systems of balance laws}
\label{cbalance}

We consider PDE systems of the form:
\begin{equation}
 w_t + F(w)_x = {S}(w)\sigma_x, \label{balance}
\end{equation}
where the unknown $w(x,t)$ takes values on an open convex set
$\mathcal{O}$ of $\mathbb{R}^{N}$;  $F$ and $S$ are regular functions from
$\mathcal{O}$ to $\mathbb{R}^{N}$; and  $\sigma(x)$ is a known function from $\mathbb{R}$ to
$\mathbb{R}$.
For an account of the existing literature on well-balanced schemes for this class of systems we refer
to the lecture notes by Bouchut \cite{Bouchut}, as well as to the recent contribution by Noelle et al. \cite{NX}, 
and Xin and Shu \cite{XS}. 

As pointed out in LeFloch \cite{LeFloch1} for the Euler equations in nozzle with discontinuous cross-section, 
such a system can be recast in the form of a nonconservative system $W_t+A(W)W_x=0$, by 
introducing 
$$
W = \left[ \begin{array}{c} w \\ \sigma \end{array} \right], \quad
A(W) = \left[
\begin{array}{c|c}
J(w) & -{S}(w)\\
\hline 0 & 0
\end{array}
\right],
$$
where $J(w)$ denotes the Jacobian matrix of $F$:
$$
{J}(w)=\frac{\partial F}{\partial w}(w).
$$
If $J$ has $N$ different real and non-vanishing eigenvalues $\lambda_1(w),\ldots,\lambda_N(w)$,
then the system is strictly hyperbolic with eigenvalues 
$$
\lambda_1(w), \dots, \lambda_N(w), 0. 
$$
Clearly, the $(N+1)$-th field is linearly degenerate and, for definiteness, we may assume that
all other fields are genuinely nonlinear.

In order to define weak solutions to this nonconservative system, a family of
paths
\[
\Phi(s; W_l, W_r) = \left[
\begin{array}{c}
\Phi_w(s; W_l, W_r)\\
\Phi_\sigma(s; W_l, W_r)
\end{array}
\right]
\]
must be chosen and it is natural to
impose the following  requirement: if $W = [ w, \sigma]^T$ is a weak
solution to the system and $\sigma $ is a constant, then $w$
must be a weak solution of the system of conservation laws:
$$
 w_t + F(w)_x = 0. 
$$
This requirement is satisfied if the family of paths
fulfills the following condition: 
\begin{itemize}

\item[(R4)] If ${W}_l$ and ${W}_r$ are such that $\sigma_l=\sigma
_r=\bar{\sigma }$, then:
\begin{equation}\label{cond_N+1}
{\Phi}_\sigma(s; {W}_l, {W}_r) = \bar{\sigma}, \qquad s \in [0,1].
\end{equation}

\end{itemize}

In fact any family of paths satisfying (R4) together with the requirements (R1) and (R3) (cf. Section \ref{schemes}), leads to the same notion of weak solution for such systems. These solutions contain two type of
discontinuities :
\begin{itemize}
\item Shock waves across which $\sigma$ is continuous that satisfy the usual Rankine-Hugoniot conditions:
\[
\xi (w^+ - w^-) = F(w^+) - F(w^-).
\]
\item Stationary contact discontinuities placed at the jumps
of $\sigma$ and connecting two states that belong to the same
integral curve of the linearly degenerate field.
\end{itemize}

For these systems, the difficulties of convergence commented above are not appreciated for the finite difference
schemes introduced in Section \ref{schemes}. Moreover, the shock waves propagating in regions where $\sigma$
is continuous are correctly captured independently of the choice of paths.  Nevertheless, in order to correctly
capture the stationary contact discontinuities related to the jumps of $\sigma$, the numerical schemes have to
be based on families of paths that satisfy at least the requirement  (R1). Interestingly, this is also the requirement
necessary to obtain well-balanced schemes.

The Lax-Friedrichs scheme presented in Section 2.2 cannot be used
as it stands for systems of balance laws  as it does not preserve
the equation
\[
\sigma_t = 0.
\]
In effect, it can be easily verified that the numerical scheme for
the variable $\sigma$ reads as follows:
\[
\sigma_i^{n+1} = \frac{\sigma_{i-1}^n + \sigma_{i+1}^n}{2}.
\]
In \cite{CPPT} the following modification of the scheme was
introduced to get rid of this difficulty: instead of
\eqref{DLFM},  $M_{i+1/2}^{n,\pm}$ are given by
$$
\aligned 
M_{i+1/2}^{n,-} & =  \widehat{A}^{n,-}_{i+1/2}\cdot (u_{i+1}^n-u_i^n), \label{DWBLFM} \\
M_{i+1/2}^{n,+} & =  \widehat{A}^{n,+}_{i+1/2}\cdot
(u_{i+1}^n-u_i^n), \label{DWBLFP}
\endaligned
$$ 
where
$$
\widehat{A}^{n,\pm}_{i+1/2} = \frac{1}{2}\left(\pm \frac{\Delta
x}{\Delta t} \widehat{I}_{i+1/2}^n+ {A}^{n}_{i+1/2}\right).
$$
Here, ${A}^{n}_{i+1/2}$ is a Roe linearization and
$$
\widehat{I}_{i+1/2}^n=\mathcal{K}_{i+1/2}^n\cdot  \widehat{Id}
\cdot (\mathcal{K}_{i+1/2}^n)^{-1}, 
$$
where $\mathcal{K}_{i+1/2}^n$ is a matrix whose columns are
eigenvectors of $A_{i+1/2}^n$ associated to $\lambda_{i+1/2,1}^n$,
\dots, $\lambda_{i+1/2,N}^n$, and $\widehat{Id}$ is the diagonal
matrix whose $j$-th coefficient is 1 if $\lambda_{i+1/2,j}^n \not=
0$, or 0 if  $\lambda_{i+1/2,j}^n = 0$. Note that if, instead,
$\widehat{I}_{i+1/2}^n$ is taken to be equal to the identity
matrix, then the Lax-Friedrichs scheme presented in Section 2.2 is
recovered.

An important particular example of systems of balance laws is the shallow water system
governing the flow of a shallow layer of inviscid homogeneous fluid through
a straight channel with a constant rectangular cross-section:
\begin{equation}
\aligned
& \frac{\partial h}{\partial t}+ \frac{\partial q}{\partial x}=0, 
\\
& \frac{\partial q}{\partial t}+ \frac{\partial}{\partial x} \left(
{\frac{q^2}{h}+  \frac{g}{2}h^2} \right) =
  gh\frac{\mbox{d} H}{\mbox{d} x}. 
\endaligned
 \label{swdv}
\end{equation}
The variable $x$ makes reference to the axis of the channel and
$t$ is time; $q(x,t)$ and $h(x,t)$ represent the mass-flow and the
thickness, respectively; $g$, the gravity; and $H(x)$, the depth
measured from a fixed level of reference.

The eigenvalues of this matrix are $\lambda_1 = u - c$, $\lambda_2 = u + c$, and $0$, 
where $c=\sqrt{gh}$. 
In this case, the equations of the integral curves of the linearly degenerate field are 
given by the equations 
\begin{equation}
q=const, \qquad h+\frac{q^2}{2gh^2}-H=const.
\label{equilibria}
\end{equation}
The resonant regime where $\lambda_1$ or $\lambda_2$ vanish are not considered in the present paper; 
for a discussion of the Riemann problem see LeFloch and Thanh \cite{LFT2}. 

%----------------------------------------------------------------------------------------------------------

\subsection{Two-layer shallow water system}
\label{twolayer}

We consider in this paragraph the system of partial
differential equations governing the one-dimensional flow of two superposed
immiscible layers of shallow water fluids  over a flat bottom topography (see \cite{CMP00a} for details) :
\begin{equation}
\aligned 
& \displaystyle{{\frac{\partial h_{1}}{\partial t}}+{\frac{\partial q_{1}}{%
\partial x}}=0,} 
\\
& \displaystyle{\frac{\partial q_{1}}{\partial t}}+{\frac{\partial \ }{%
\partial x}}\left( {\frac{q_{1}^{2}}{h_{1}}}+{\frac{g}{2}}h_{1}^{2}\right)
=-gh_{1}{\frac{\partial h_{2}}{\partial x}}, 
\\
& \displaystyle{{\frac{\partial h_{2}}{\partial t}}+{\frac{\partial q_{2}}{%
\partial x}}=0,} 
\\
& \displaystyle{{\frac{\partial q_{2}}{\partial t}}+{\frac{\partial \ }{%
\partial x}}\left( {\frac{q_{2}^{2}}{h_{2}}}+{\frac{g}{2}}h_{2}^{2}\right)
=-{\frac{\rho_{1}}{\rho_{2}}}gh_{2}{\frac{\partial
h_{1}}{\partial x}}}.
\endaligned
   \label{MSW}
\end{equation}
In these equations, the index 1 refers to the upper layer and the index 2
to the lower one. The fluid is assumed to occupy a straight channel with
constant rectangular cross-section and constant width. The  coordinate
$x$ refers to the axis
of the channel, $t$ denotes the time variable, and $g$ is the gravity.  Each layer is assumed to have a constant
density $\rho_{i}$, $i=1,2$ ($\rho_{1}<\rho_{2}$), while 
the unknowns $q_{i}(x,t)$ and $h_{i}(x,t)$ represent respectively the mass-flow and the thickness of the $i$-th layer at
the section of coordinate $x$ at time $t$.

System (\ref{MSW}) can be written in the form (\ref{NONC}), say 
$w_t + A(w) \, w_x = 0$, $w=w(t,x) \in \RN$,
with $N = 4$ and 
$$
w= \left[
\begin{array}{c}
h_1\\
q_1\\
h_2\\
q_2\\
\end{array}
\right], 
\qquad 
A(w) = \left[
\begin{array}{cccc}
0 & 1 & 0 & 0 \\
-u_1^2 + c_1^2 & 2u_1 & c_1^2 & 0 \\
0 & 0 & 0 & 1 \\
rc_2^2 & 0 & -u_2^2 + c_2^2 & 2u_2\\
\end{array}
\right],
$$
where $u_i =  {q_{i}}/{h_{i}}$ represents the averaged velocity of
the $i$-th layer, $c_i = \sqrt{gh_i}$, $i = 1, 2$, and $r=\frac{\rho_1}{\rho_2} $.
The characteristic equation of the system is 
$$
\bigl( \lambda^2 - 2u_1\lambda + u_1^2 - g h_1 \bigr)\bigl( \lambda^2 -
2u_2\lambda + u_2^2 - g h_2 \bigr) = r g^2 h_1 h_2. 
$$
Observe that, when $r=0$, the eigenvalues are those corresponding to
each layer separately.   In this situation, the coupling terms do not affect the nature of the
system in any essential manner.

In the case $r\cong 1$ (which is the situation arising in many geophysical flows)
a first-order approximation of the eigenvalues was
given in \cite{SS53}:
\begin{eqnarray}
\lambda_{\hbox{\scriptsize{ext}}}^{\pm} & \cong & {\frac{u_{1}h_{1}+u_{2}h_{2}}{h_1+h_2}} \pm
         \bigl(g(h_{1}+h_{2})\bigr)^{1/2},  \label{eq_lambda_ext} \\
\lambda_{\hbox{\scriptsize{int}}}^{\pm} & \cong & {\frac{u_{1}h_{2}+u_{2}h_{1}}{h_1+h_2}}\pm
         \left( g^{\prime }{\frac{h_{1}h_{2}}{(h_{1}+h_{2})}}
         \Bigl( 1-{\frac{(u_{1}-u_{2})^{2}}{g^{\prime}(h_{1}+h_{2})}}
         \Bigr)\right)^{1/2}.
         \label{eq_lambda_int}
\end{eqnarray}
In the former expression, $g^{\prime }=(1-r)g$  is the \textit{reduced gravity}. 

It is not easy to check the genuinely nonlinear character of the 4 characteristic fields,
as the eigenvalues and
eigenvectors can not be written  explicitely in a simple manner. Nevertheless, this fact is
easily proved in the case $r=0$ as, in this case, the system reduces to two decoupled
shallow water systems. As a consequence, using a continuity argument, this is also true
at least for small values
of $r$.

From equation (\ref{eq_lambda_int}) we can observe that the internal
eigenvalues may become complex.
This situation occurs when they satisfy, approximately, the following
inequality:
$$
{\frac{(u_{1}-u_{2})^{2}}{g^{\prime }(h_{1}+h_{2})}}>1.
$$
In this case, the system loses its hyperbolic character.
These situations are related with the appearance of shear
instabilities that may lead,  in real flows, to
intense mixing of the two layers.
While, in practice, this mixture partially dissipates the energy, in
numerical experiments these interface disturbances grow and overwhelm
the solution.
Clearly, we cannot expect to simulate these phenomena with a
two-immiscible-layer model.
Therefore, the above inequality in fact gives the range of
validity of a model based on the equations (\ref{MSW}), if viscosity
effects are neglected. In this work only the case where the matrix ${A}({w})$ has
 real eigenvalues is considered, i.e. the system is supposed to
be strictly  hyperbolic.

The jump conditions \eqref{RH} related to the choice of a family of paths
$$
\Phi(s; w_l, w_r) =
\left[
\begin{array}{c}
\Phi_{h_1}(s; w_l, w_r) \\
\Phi_{q_1}(s; w_l, w_r) \\
\Phi_{h_2}(s; w_l, w_r)\\
\Phi_{q_2}(s; w_l, w_r)
\end{array}
\right]
$$
read 
\begin{equation}\label{RH-2C-1}
\aligned 
\xi \displaystyle{(h_1^r - h_1^l)} 
= & \displaystyle{q_1^r  - q_1^l},
 \\
\xi \displaystyle{(q_1^r - q_1^l)} 
= & \displaystyle{\frac{(q_1^r)^2}{h_1^r}  - \frac{(q_1^l)^2}{h_1^l}
 + \frac{g}{2}(h_1^r)^2 -  \frac{g}{2}(h_1^l)^2} 
 \\
& \displaystyle{+ g \int_0^1 \Phi_{h_1}(s; w_l, w_r)\frac{\partial  \Phi_{h_2} }{\partial s}(s; w_l, w_r)\, ds},
\\
\xi \displaystyle{(h_2^r - h_2^l)} 
= &  \displaystyle{q_2^r - q_2^l},
 \\
\xi \displaystyle{(q_2^r - q_2^l)} 
= & \displaystyle{\frac{(q_2^r)^2}{h_2^r} - \frac{(q_2^l)^2}{h_2^l}
    + \frac{g}{2}(h_2^r)^2 -  \frac{g}{2}(h_2^l)^2}\\
 & \displaystyle{+gr \int_0^1 \Phi_{h_2}(s; w_l, w_r)
\frac{\partial  \Phi_{h_1}}{ \partial s}(s; w_l, w_r)\,ds.}
\endaligned 
\end{equation}
Observe that these conditions are independent of the choice of $\Phi_{q_i}(s: w_l, w_r)$, $i=1,2$.

%*****************************************************************************************

\begin{center}
\begin{figure}[h!ptb]
\centering
\includegraphics[width=8.5cm,height=8.0cm]{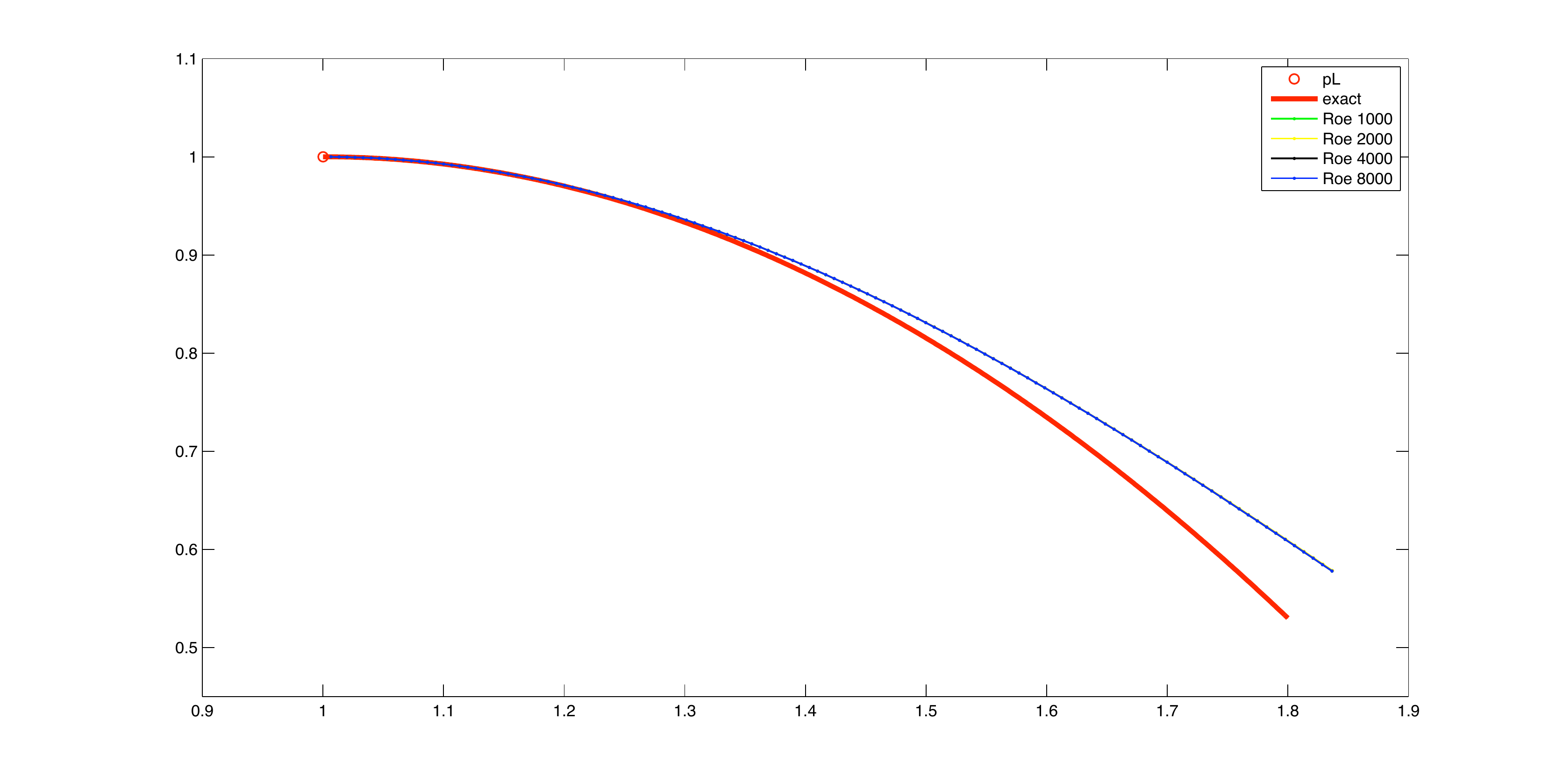}
\caption{Test case \ref{testsimple}. Hugoniot curves: exact (continuous red line) and numerical (line with dots) obtained with Roe scheme}
\label{test_s1_fig}
\end{figure}
\end{center}

\section{Numerical experiments}

\subsection{A simplified system}
\label{testsimple}

We will  demonstrate here that, in general, a difference scheme for a general nonconservative hyperbolic
system \eqref{NONC} does not converge to the exact solution $u$, that is, with the notation already introduced in previous
sections we claim that
\be
\label{no}
v = \lim_{\Delta \to 0} u^\Delta \neq u.
\ee
This fact was first observed for nonconservative schemes for scalar equations in Hou and LeFloch \cite{HouLeFloch} and, in the
context of nonclassical shocks generated by diffusion and dispersion, in Hayes and LeFloch \cite{HayesLeFloch,LM}. 
Here, we observe \eqref{no} for finite difference approximations of nonconservative systems.

\begin{center}
\begin{figure}[h!ptb]
\centering
\includegraphics[width=8.5cm,height=8.0cm]{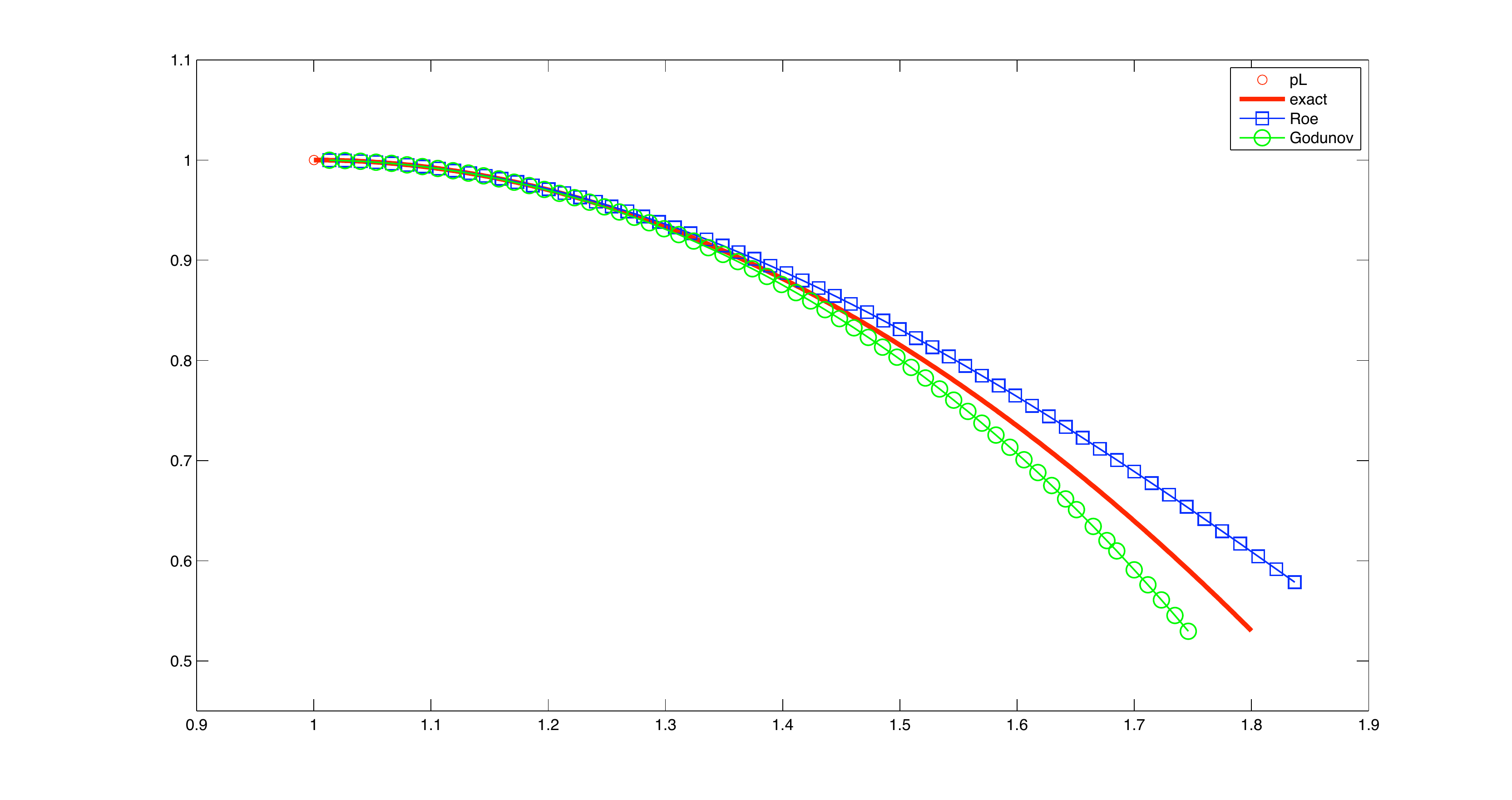}
\caption{Test case \ref{testsimple}. Hugoniot curves: exact (continuous red line), Roe (blue line with squares)  and Godunov (green line with circles) }
\label{test_s2_fig}
\end{figure}
\end{center}

We begin with the convergence of Roe and Godunov methods for \eqref{simple}.
The Roe method considered here is consistent with the family of paths given by \eqref{camino_simple} for every 
pair of states. The corresponding Roe matrix is as follows:
$$
A_{i+1/2}^n= \left[
\begin{array}{cc}
0 & 1 \\
-(u_{i+1/2}^n)^2 +q_i^n h_{i+1/2}^n & 2u^n_{i+1/2}
\end{array}
\right] ,
$$
where
$$
u_{i+1/2}^n=\frac{\sqrt{h_i^n} u_i^n+\sqrt{h_{i+1}^n}u_{i+1}^n}{\sqrt{h_i^n}+\sqrt{h_{i+1}^n}},
\quad h_{i+1/2}^n = \frac{1}{2}(h_i^n + h_{i+1}^n).
$$

On the other hand, the Godunov method considered here is based on the family of paths  described in 
Section \ref{csimple}: it satisfies (R2)-(R3) and coincides with \eqref{camino_simple} for pair of states that can be
linked by a shock.

We consider the  Rankine-Hugoniot curve composed by the states $w_r$ that can connected with $w_l=[1 ;1]^T$ by a 1-shock,
which is given by:
\begin{equation} \label{Hugoniot-1}
q_r=h_r\left(1-\sqrt{\frac{h_r+1}{2h_r}}(h_r-1)\right).
\end{equation}
Figure \ref{test_s1_fig} shows a plot of the curve \eqref{Hugoniot-1} in the plane $h-q$ (continuous red line).

\begin{center}
\begin{figure}[h!ptb]
\subfigure[h at $t=0.5$]{\
\begin{minipage}[b]{0.48\textwidth}
\centering
\includegraphics[width=6.5cm,height=6.0cm]{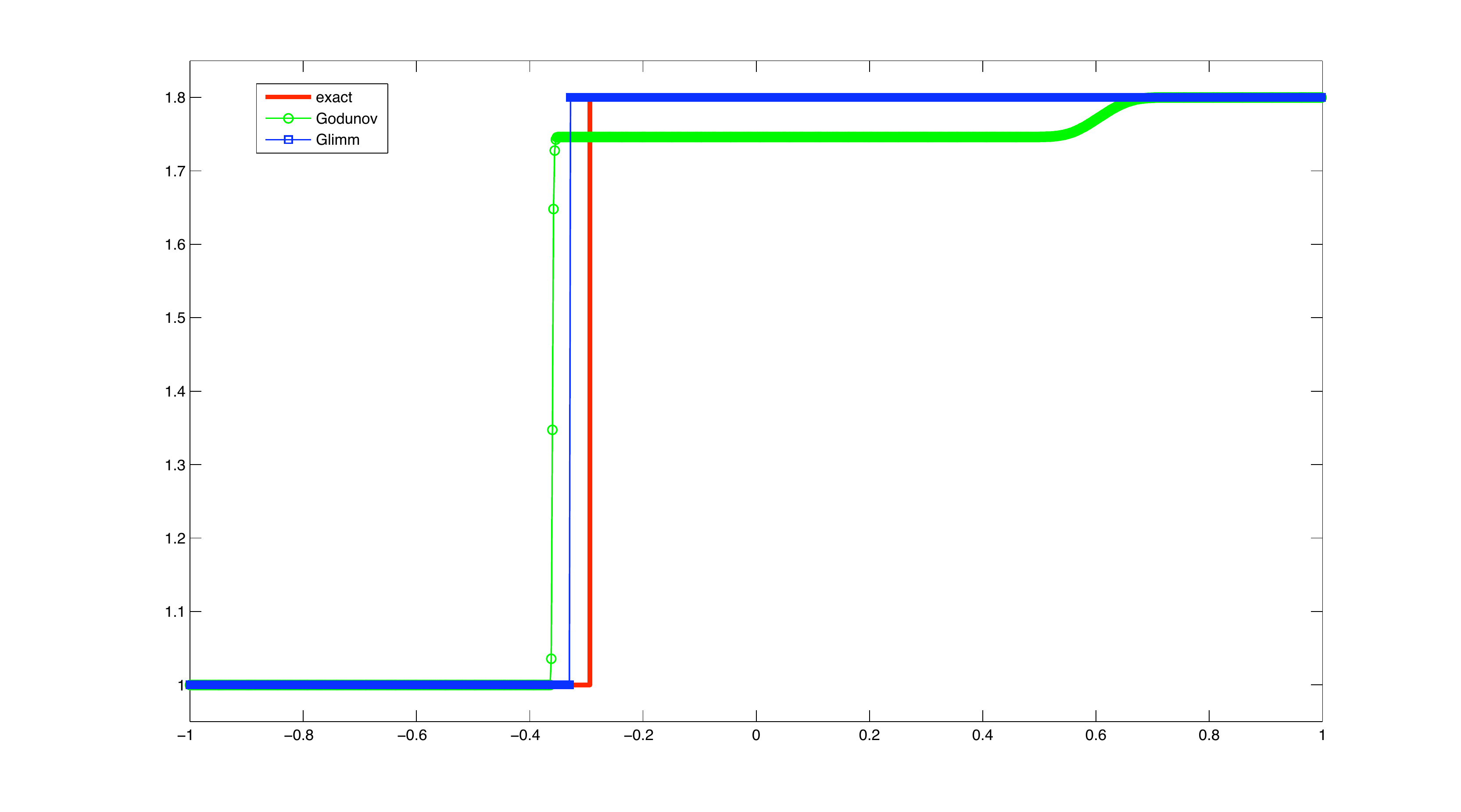}
\end{minipage}}
\subfigure[h (zoom) at $t=0.5$.]{\
\begin{minipage}[b]{0.48\textwidth}
\centering
\includegraphics[width=6.5cm,height=6.0cm]{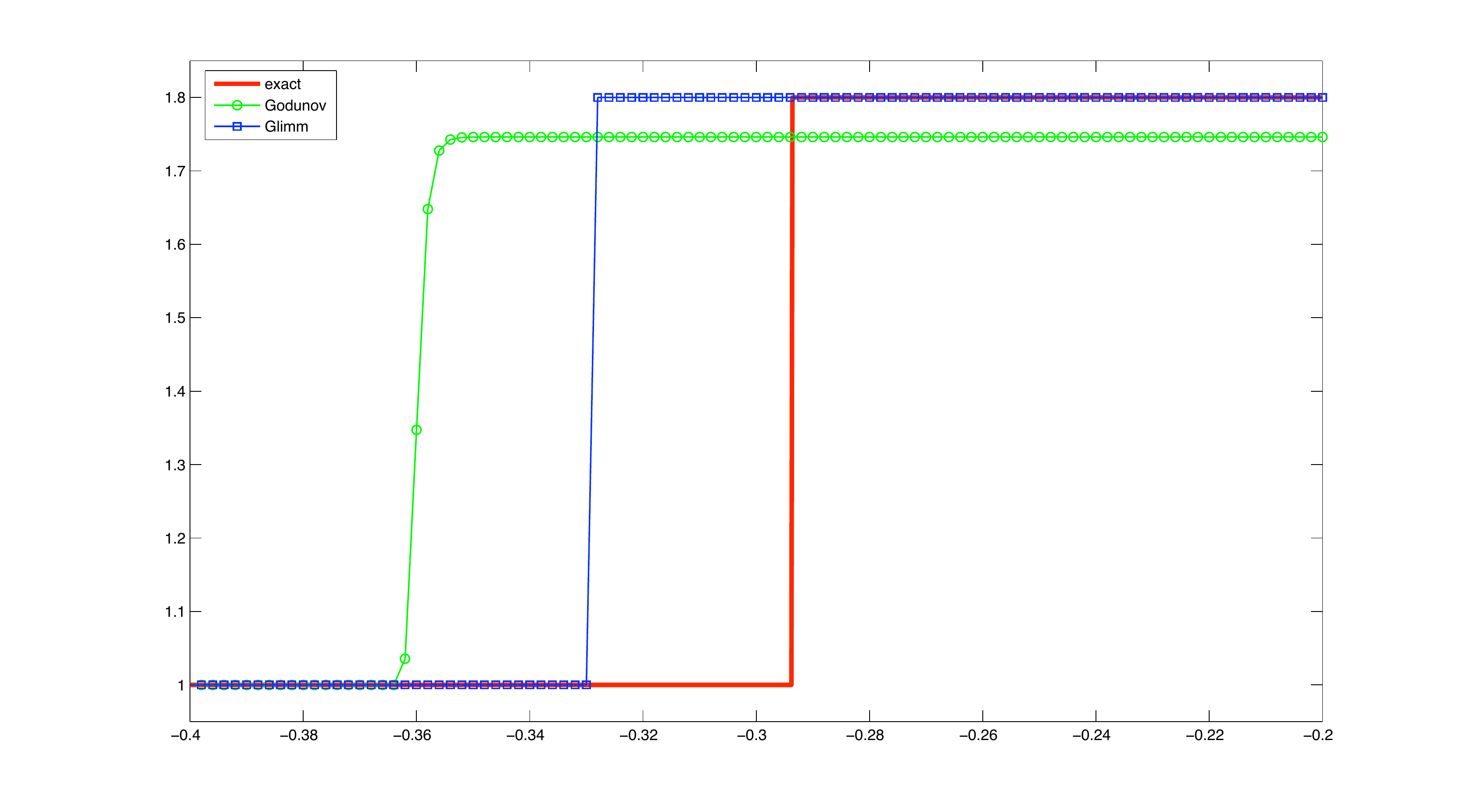}
\end{minipage}}
\caption{Test case \ref{testsimple}. Solution of Riemann problem ($h$) \eqref{rmp} at time $t=0.5$: Exact (continuous red line), Godunov (green line with circles) and  Glimm (blue line with squares) }
\label{test_s3_fig}
\end{figure}
\end{center}

\begin{center}
\begin{figure}[h!ptb]
\subfigure[q at $t=0.5$]{\
\begin{minipage}[b]{0.48\textwidth}
\centering
\includegraphics[width=6.5cm,height=6.0cm]{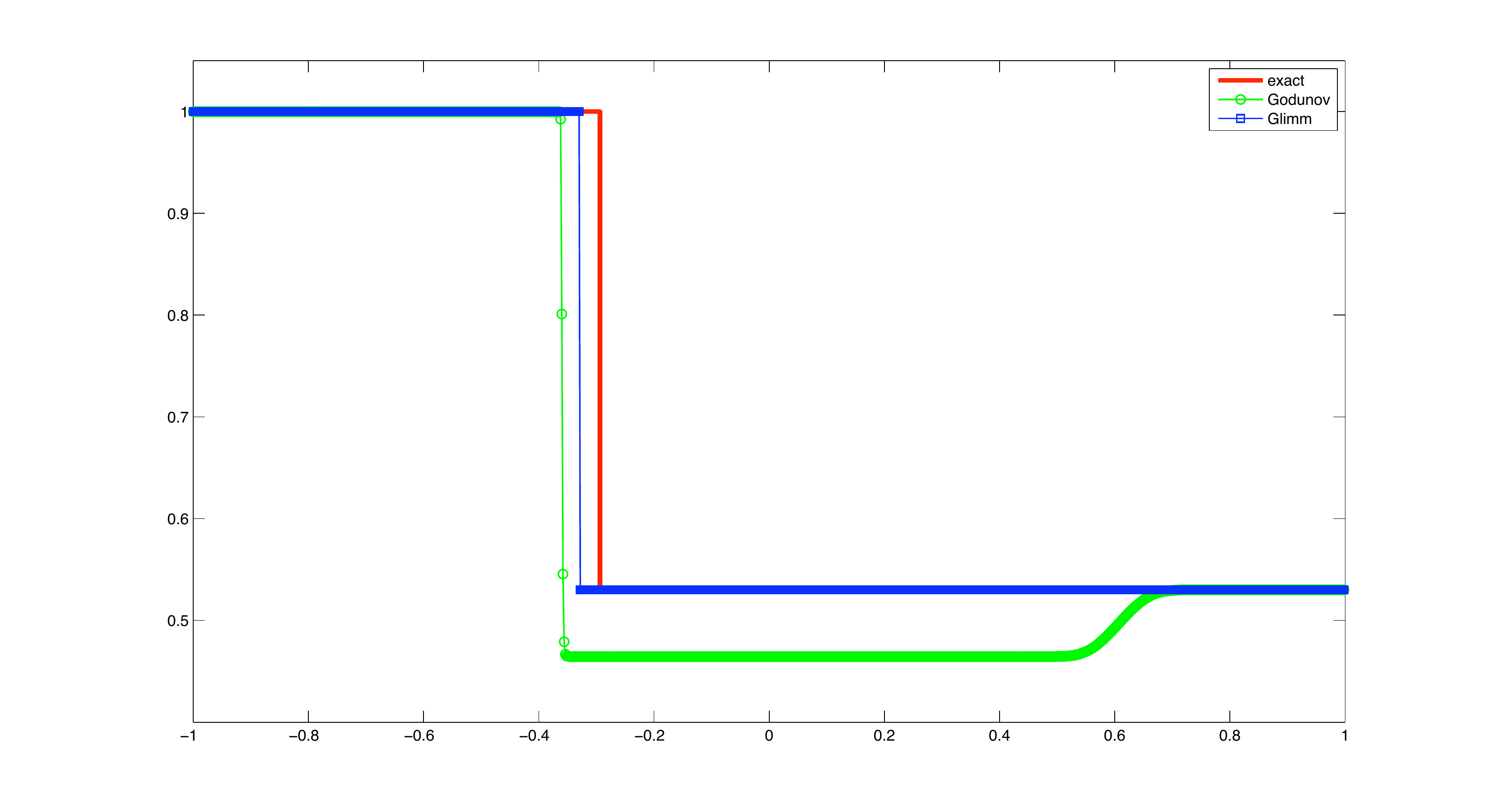}
\end{minipage}}
\subfigure[q (zoom) at $t=0.5$.]{\
\begin{minipage}[b]{0.48\textwidth}
\centering
\includegraphics[width=6.5cm,height=6.0cm]{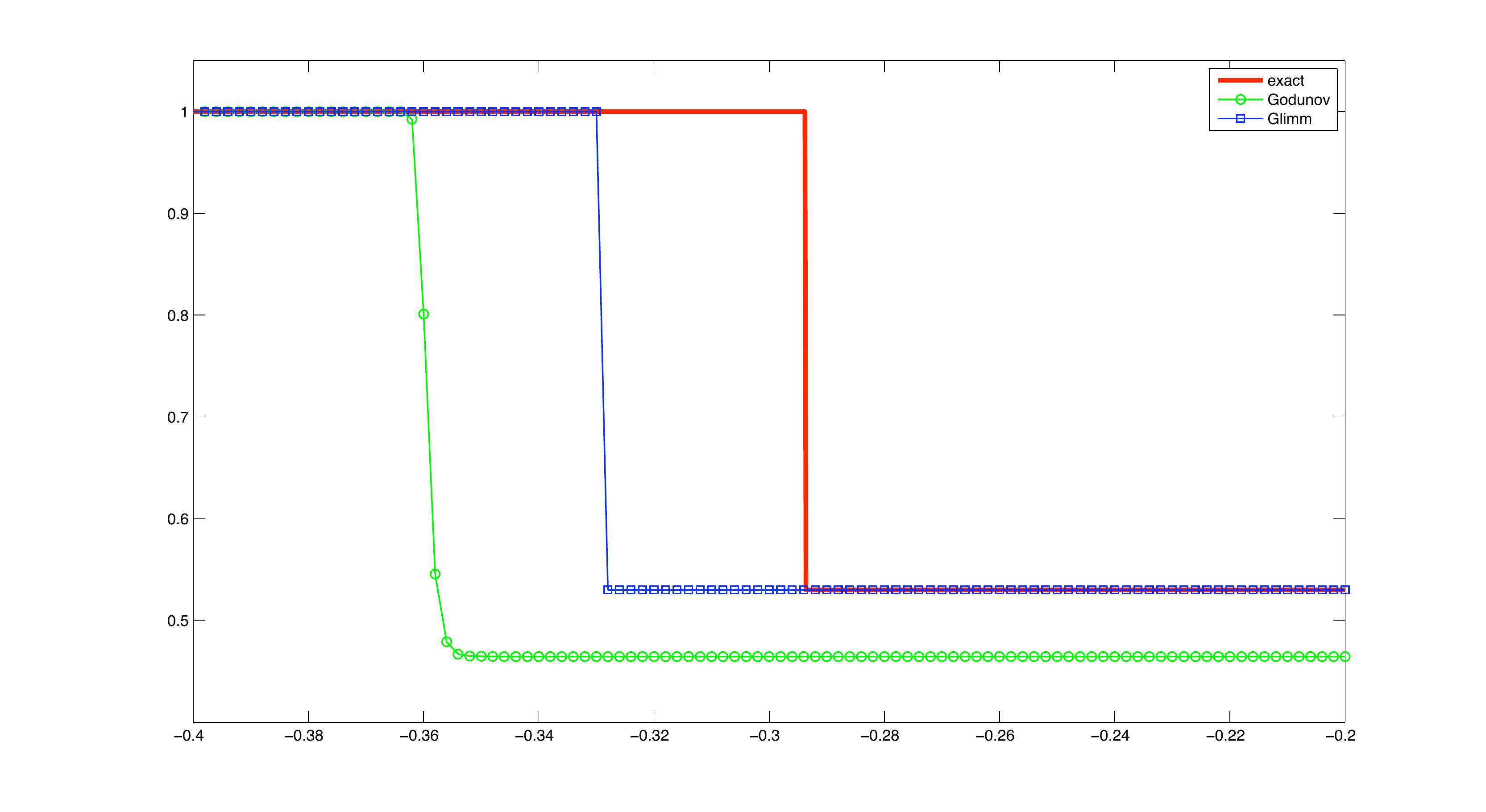}
\end{minipage}}
\caption{Test case \ref{testsimple}. Solution of Riemann problem ($q$) \eqref{rmp} at time $t=0.5$: Exact (continuous red line), Godunov (green line with circles) and  Glimm (blue line with squares) }
\label{test_s4_fig}
\end{figure}
\end{center}

Next,  Roe method is used to solve numerically a family of Riemann problems in which the left state is $w_l$ and the right state $w_r$ runs on the Hugoniot curve \eqref{Hugoniot-1}. The speed of propagation and the limit states of the 1-shock
are computed in the numerical solution by using the first divided difference as a smooth indicator. This computation has been performed using  four meshes with decreasing mesh step ($\Delta x=0.002$, $\Delta x=0.001$, $\Delta x=0.0005$ and $\Delta x=0.00025$). The CFL parameter is set to $0.9$. The numerical Hugoniot curves obtained in this way (line with dots) are compared with the exact one (continuous red line) in Figure \ref{test_s1_fig}.  It can be observed that the numerical 
Hugoniot curves converge, but the limit is not  the exact one.  

The same behavior is observed for Godunov method: the numerical Hugoniot-curves converge but the limit is not the exact one:
in Figure \ref{test_s2_fig}  the exact Hugoniot curve (continuous red line) is compared with those computed with Roe (blue line with squares) and Godunov (green line with circles)  methods  with $\Delta x=0.001$. The CLF parameter is set to $0.9$ for Roe and $0.5$ for Godunov. This choice ensures that the Godunov method corresponds to advance in time by exactly solving the Riemann problems and taking the averages of the solutions at the cells.

Finally, we compare Godunov and Glimm methods. We consider  a Riemann problem with initial conditions
\begin{equation} \label{rmp}
w(x,0)=\left\{
\begin{array}{lr}
w_l = \left( \begin{array}{c} 1 \\ 1 \end{array} \right),Ê&  x< 0; \\
\\
w_r = \left( \begin{array}{c} 1.8 \\ q_r \end{array} \right),Ê\qquad 
& x>0;
\end{array}
\right.
\end{equation}
where $q_r$ is given by \eqref{Hugoniot-1} for $h_r = 1.8$ ($q_r \cong 0.530039370688997$).
The exact solution  consists thus of a 1-shock linking the states.  In  Figures \ref{test_s3_fig} and \ref{test_s4_fig} the exact solution at time $t=0.5$ is compared with numerical solutions obtained with Godunov (green line with circles) and Glimm  (blue line with squares) methods with $\Delta x=0.001$ and CFL=$0.5$. Note how Godunov method introduces a 2-rarefaction in the computed solution.  

\begin{center}
\begin{figure}[h!ptb]
\centering
\includegraphics[width=8.5cm,height=8.0cm]{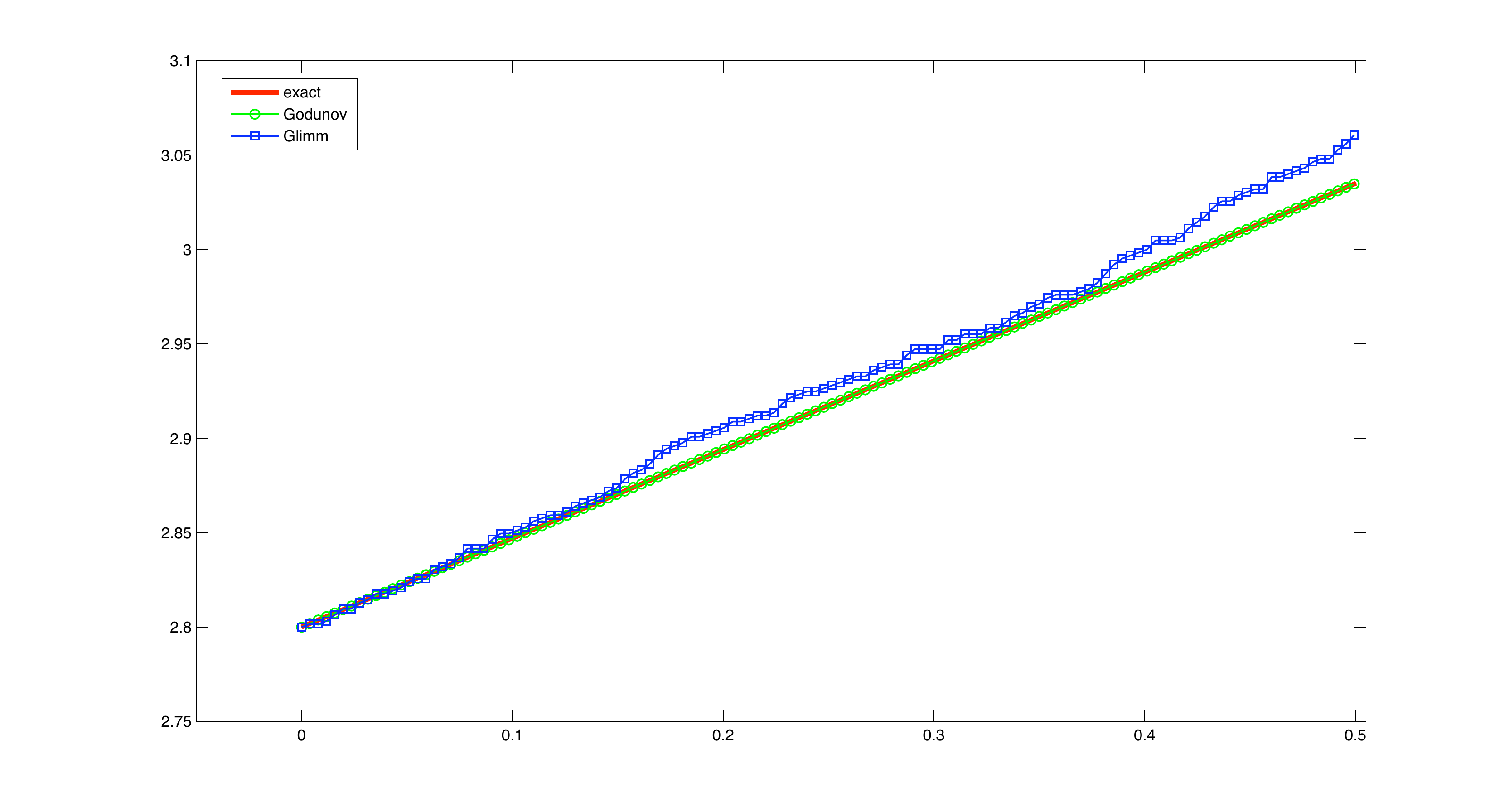}
\caption{Test case \ref{testsimple}. Evolution of the \lq mass\rq: exact (continuous red line), Glimm (blue line with squares)  and Godunov (green line with circles) }
\label{test_s5_fig}
\end{figure}
\end{center}

Observe that the first equation of  \eqref{simple} is a conservation law. According to this equation, 
for $A> 0$ large enough, the exact solution of the Riemann
problem has to satisfy the following conservation property:
$$
\int_{-A}^A h(x,t) \,dx = \int_{-A}^A h(x,0)\, dx + t (1 - q_r).
$$
 In Figure  \ref{test_s5_fig} we investigate the conservation property of both
Godunov and Glimm methods: we fix $A$ and compare the exact value of 
the integral of $h$ at time $t_n$ with its numerical approximations. Note that Godunov method satisfies 
the conservation property exactly.

%-----------------------------------------------------------------------------------------------------

\subsection{Shallow water system}
\label{sw}

In this section we consider the discretization of the shallow water system by means of Roe and the modified Lax-Friedrichs schemes applied to the formulation \eqref{NONC} of the problem. In the first test, we check  that, for continuous bottom
functions $H$, the  shock waves are correctly captured for the schemes even if the numerical schemes are
formally consistent with a simple family of paths.  In the second one,  stationary contact discontinuities placed at the jumps
of $H$ are considered: we check that they are correctly captured only when the family of paths satisfies the property
(R1).

\subsubsection{Dam-break problem over a non-flat bottom topography}
\label{test2.1}

The axis of the channel is the interval
$[0, 10]$ and the bottom topography is given by the function
$\displaystyle
H(x)=
1-0.5e^{-(x-5)^2}$.
The initial condition is $q=0$ and
$$
h(x)=\left\{
\begin{array}{ll}
\displaystyle
H(x),  \qquad & \qquad  \text{$ x\geq4$};\\ \\
H(x)+0.5, \qquad  & \qquad \text{$ x<4$}.
\end{array}
\right.
$$
The final time is $t=0.6$. Free boundary conditions are considered. The CFL parameter is set to $0.9$.

\begin{center}
\begin{figure}[h!ptb]
\subfigure[Bottom topography and free surface at t=0.6.]{\
\begin{minipage}[b]{0.48\textwidth}
\centering
\includegraphics[width=6.5cm,height=6.0cm]{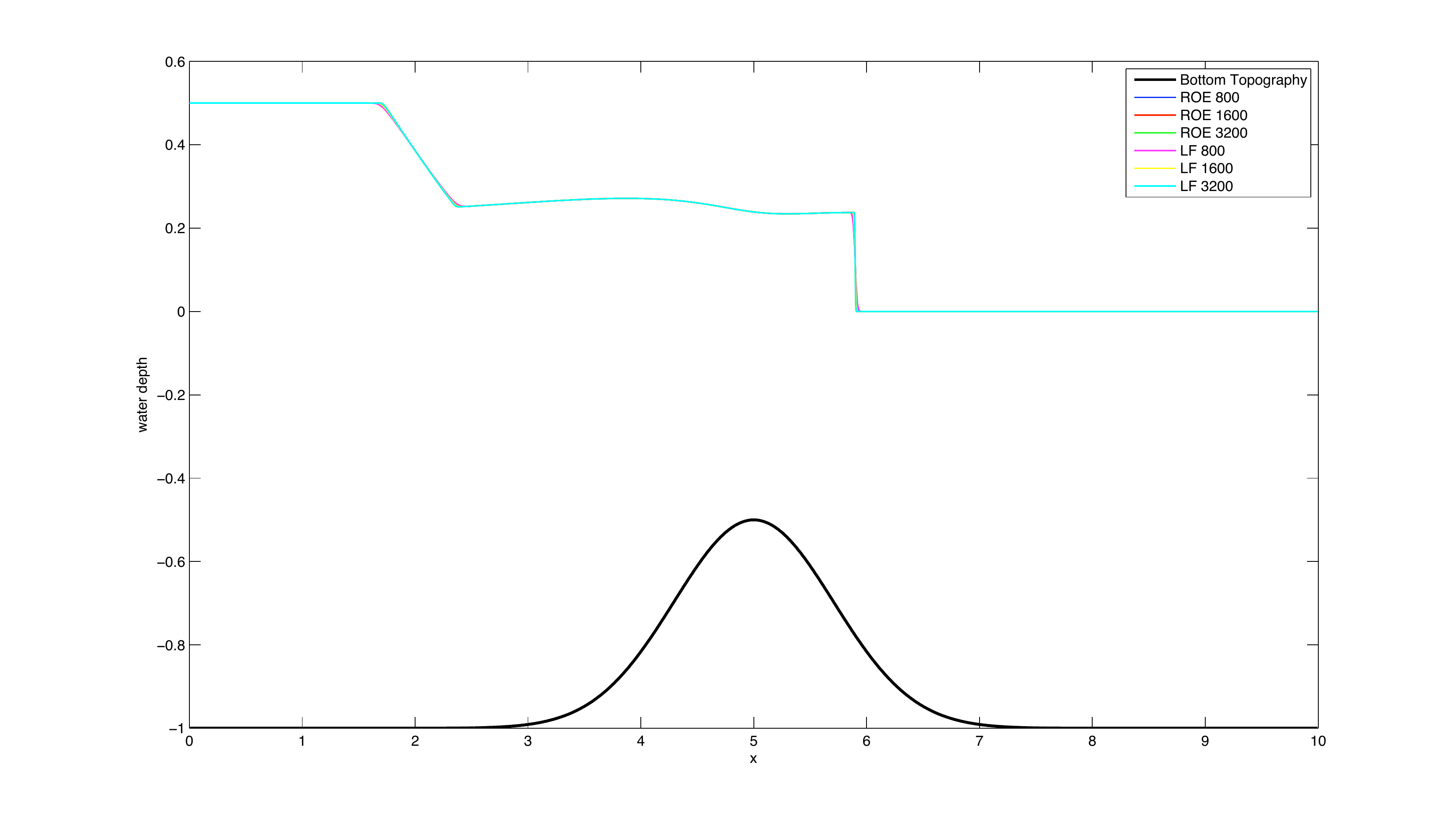}
\end{minipage}}
\subfigure[Free surface (zoom) at t=0.6.]{\
\begin{minipage}[b]{0.48\textwidth}
\centering
\includegraphics[width=6.5cm,height=6.0cm]{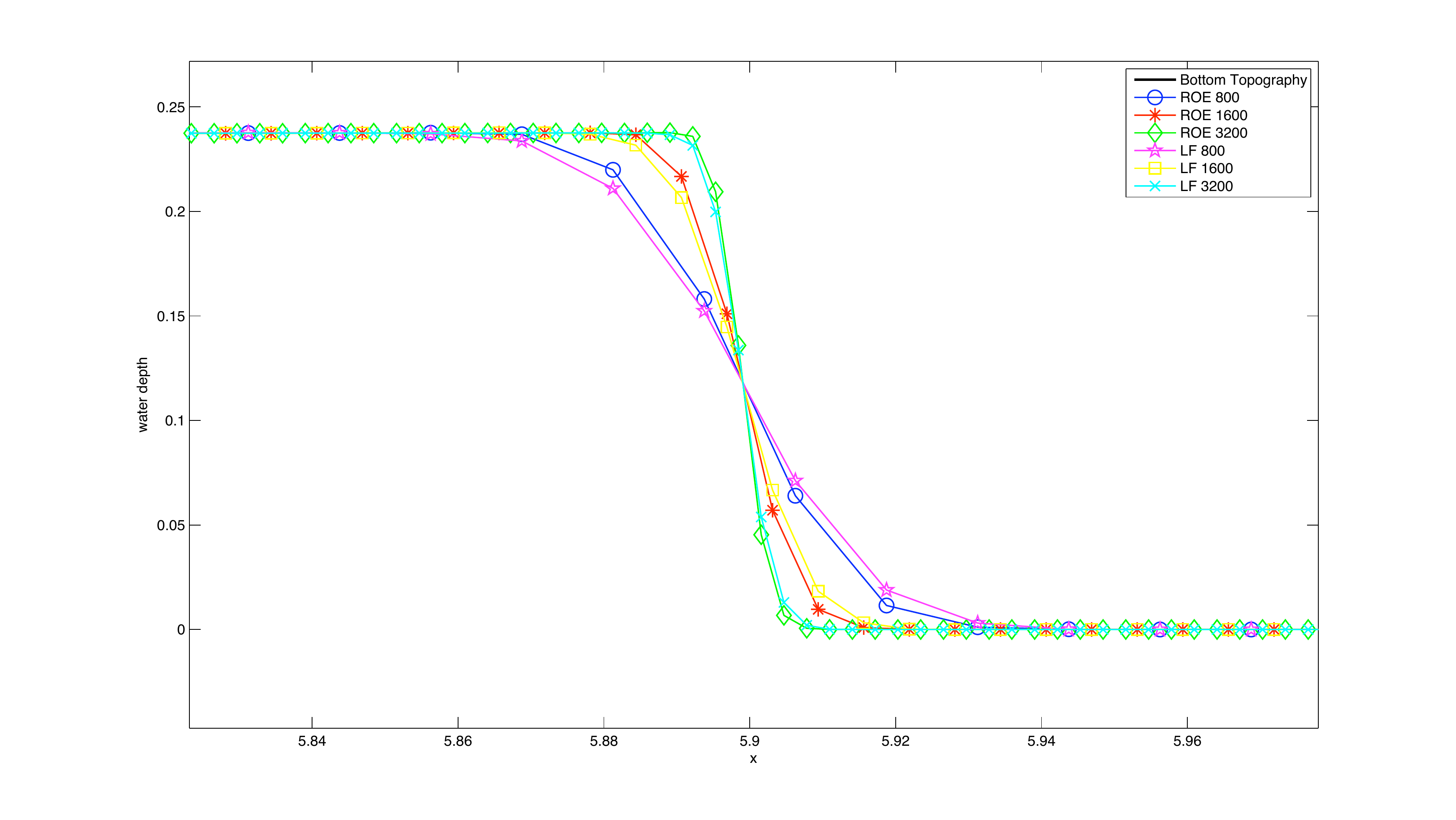}
\end{minipage}}
\caption{Test case \ref{test2.1}. Dam-break problem: bottom topography and free surface at t=0.6.}
\label{test21_fig}
\end{figure}
\end{center}

We consider a Roe scheme and a modified Lax-Friedrichs scheme which are consistent with the family of
segments.  Figure \ref{test21_fig} shows the bottom topography and the free surface computed for both schemes using three meshes with increasing number of cells (800, 1600 and 3200 cells respectively). Both schemes converge to the same solution. Moreover, the speed of propagation $\xi$ and the limit states $w^-$ and $w^+$  of the shock in the numerical solutions have
been computed for both schemes by using a fine mesh of 32000 cells and the first divided difference  as a smooth indicator. 
The value of the residual $|\xi (w^+ - w^-) - F(w^+) + F(w^-)|$ obtained for the well-balanced Lax-Friedrichs scheme 
is $0.008$  and  for Roe scheme, $0.006$.

\subsubsection{Stationary contact discontinuities}
\label{test2.2}

In this test we study the approximation of stationary contact discontinuities.
Following
the discussion in
Section \ref{cbalance}, such a discontinuity has to connect two states belonging to a same curve of
the family \eqref{equilibria}.  As the family of segments only satisfies (R1) for the particular case
of the curves
$$
q=0, \ h-H=const,
$$
the Roe and the modified Lax-Friedrichs schemes based on the family of segments are only expected to
capture correctly stationary contact discontinuities corresponding to water at rest over a discontinuous bottom.
To check this in practice,
we consider a channel whose axis is the interval $[-5,5]$ and
whose bottom is given by the function
\[H(x)=\left\{\begin{array}{cll}
0,\hspace{0.2cm}& &  x<0,
\\
1,\hspace{0.2 cm}&  & x>=0.
\end{array}\right.
\]
We consider the initial condition:
\[w_l=\left[\begin{array}{c} h_l\\ q \end{array}\right]=\left[\begin{array}{c} 1\\ \sqrt{4g} \end{array}\right], 
\qquad \quad 
w_r=\left[\begin{array}{c} h_r\\ \sqrt{4g} \end{array}\right],
\]
where $h_r$ has been calculated so that both states
belong to the same integral curve \eqref{equilibria} ($h_r \cong 0.7892441190408083$).
The exact solution of this Riemann problem is thus a stationary contact discontinuity.

\begin{center}
\begin{figure}[h!ptb]
\subfigure[Bottom topography and free surface (stationary solution).]{\
\begin{minipage}[b]{0.48\textwidth}
\centering
\includegraphics[width=6.5cm,height=6.0cm]{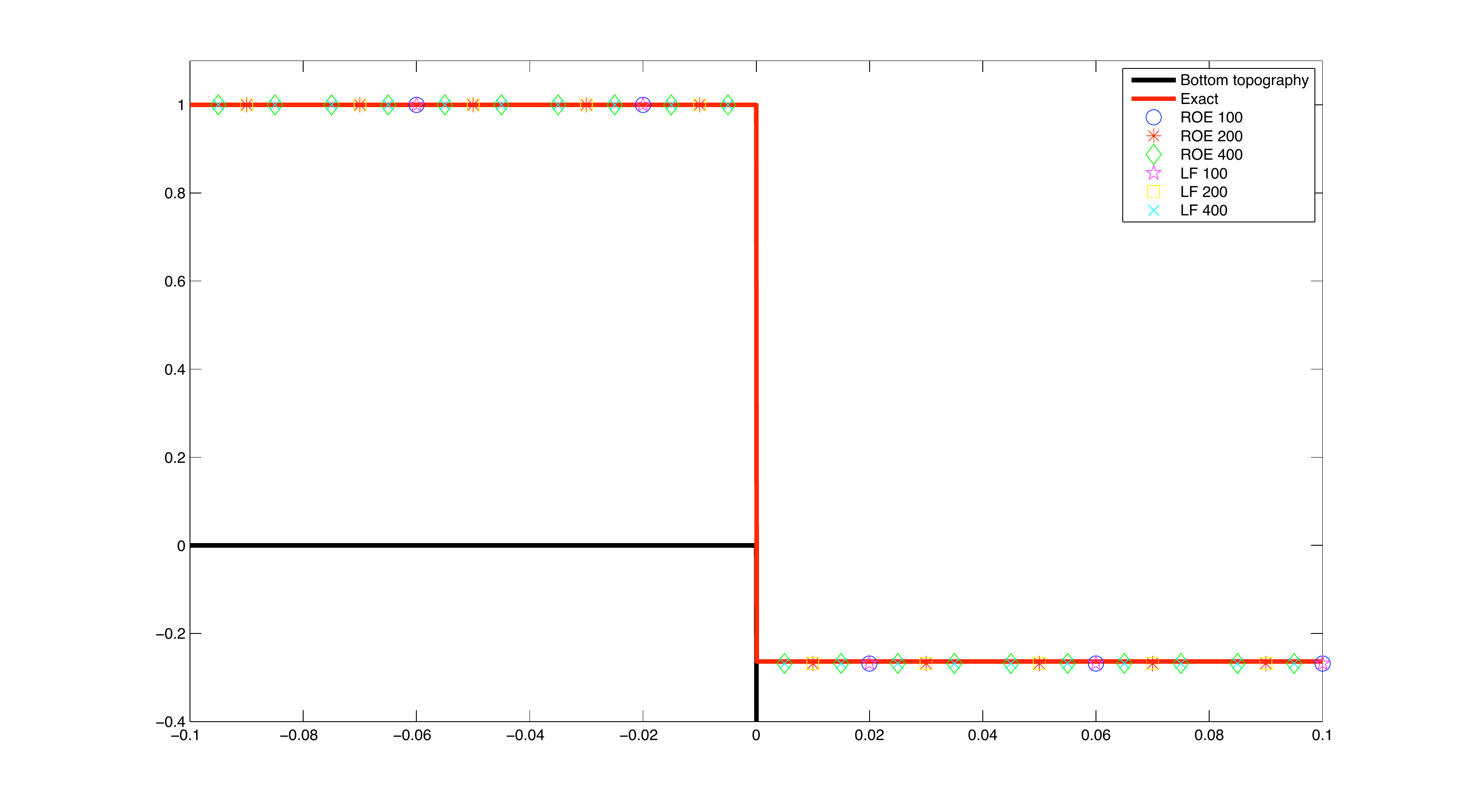}
\end{minipage}}
\subfigure[Bottom topography and free surface (stationary solution): zoom]{\
\begin{minipage}[b]{0.48\textwidth}
\centering
\includegraphics[width=6.5cm,height=6.0cm]{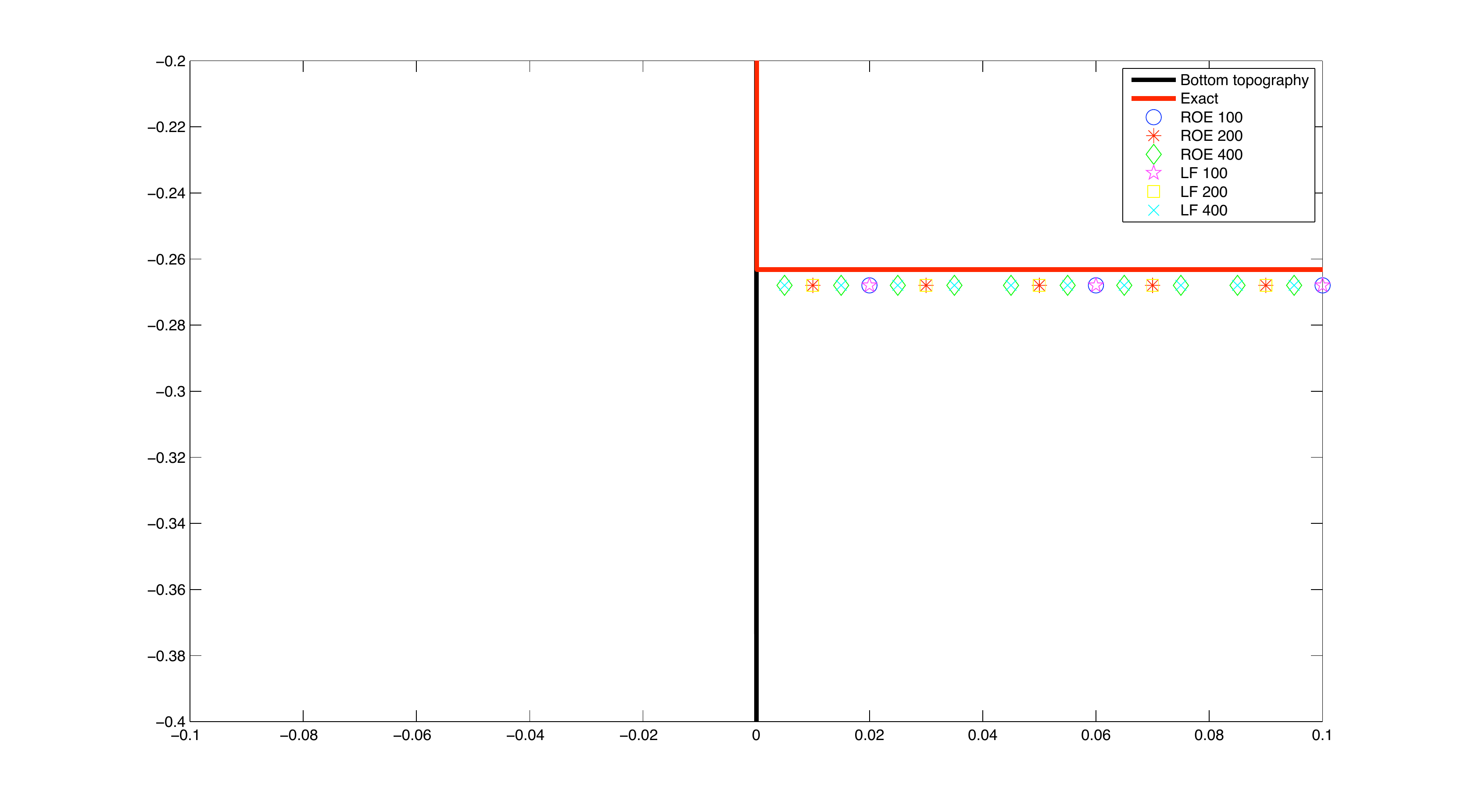}
\end{minipage}}
\caption{Test case \ref{test2.2}. Stationary contact discontinuity: Comparison between the modified Lax-Friedrichs and Roe schemes and the exact solution (free surface).}
\label{test22_fig1}
\end{figure}
\end{center}

\begin{center}
\begin{figure}[h!ptb]
\centering
\includegraphics[width=6.5cm,height=6.0cm]{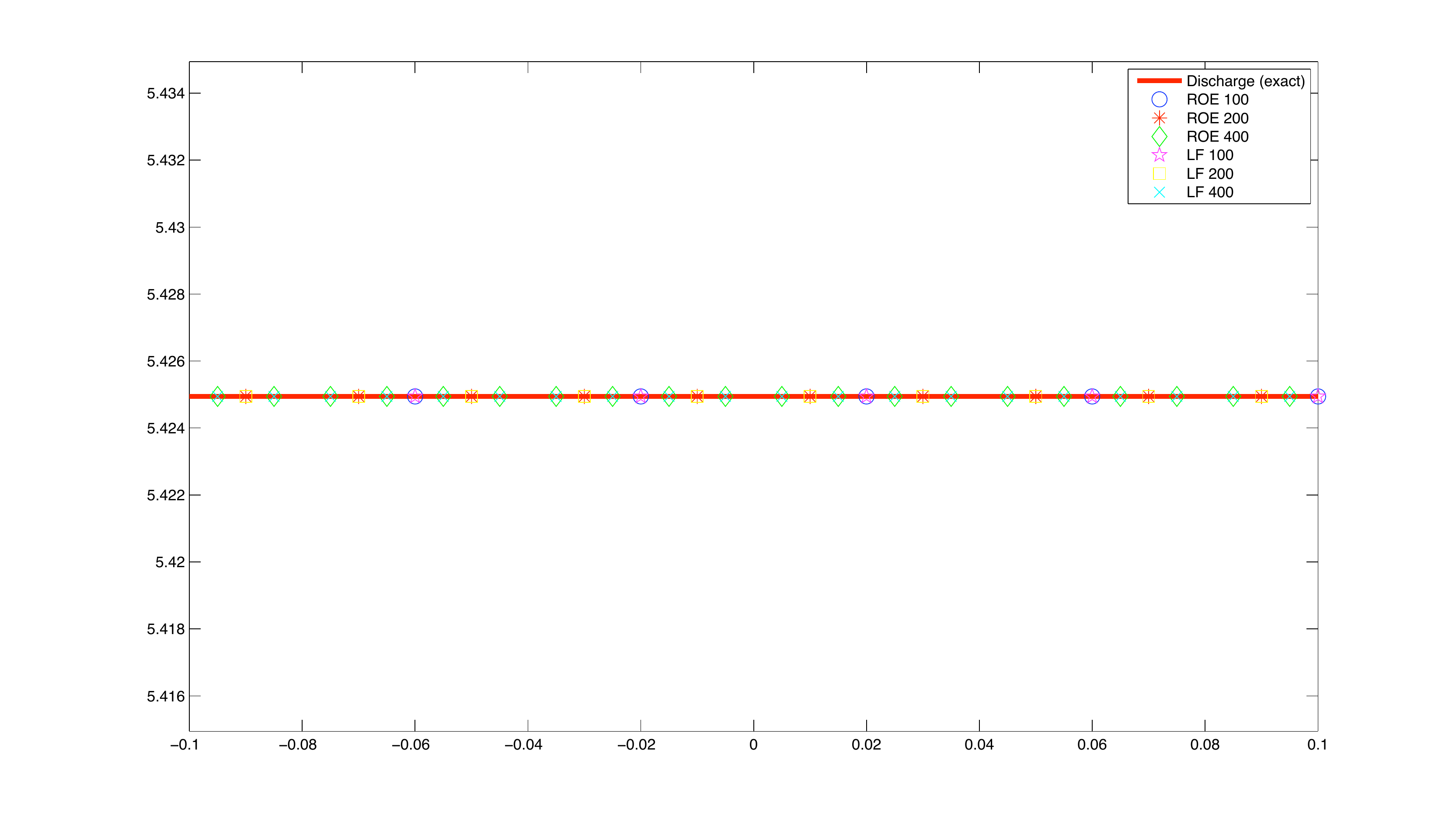}
\caption{Test case \ref{test2.2}. Stationary contact discontinuity: Comparison between the well-balanced Lax-Friedrichs and Roe schemes and the exact solution (discharge).}
\label{test22_fig2}
\end{figure}
\end{center}

We have applied both schemes to this Riemann problem. As boundary condition, the state $w_l$ is imposed upstream  and free boundary conditions  downstream. The CLF parameter is set to 0.9.
Figures \ref{test22_fig1} and \ref{test22_fig2} show the stationary solutions obtained with both schemes using three meshes with increasing number of cells (100, 200 and 400 cells respectively). As expected, the numerical solutions do not converge to the
exact solution. Note that both schemes converge to the same discontinuous function.

Let us check
 what happens if the numerical schemes are based on a family of paths satisfying (R1) for every integral curve
of the linearly degenerate field. We consider the family of paths constructed in   \cite{CPP07}  in order to design a Roe
scheme which is well-balanced for every smooth stationary solution. The path connecting two states
$W_l$ and $W_r$ consists of an arc of one of the curves of the family \eqref{equilibria} and a segment lying on a
plane $H = const$: see  \cite{CPP07}  for details. Once the corresponding Roe matrix has been calculated, the construction of the modified Lax-Friedrichs scheme  is straightforward.

\begin{center}
\begin{figure}[h!ptb]
\subfigure[Bottom topography and free surface (stationary solution).]{\
\begin{minipage}[b]{0.48\textwidth}
\centering
\includegraphics[width=6.5cm,height=6.0cm]{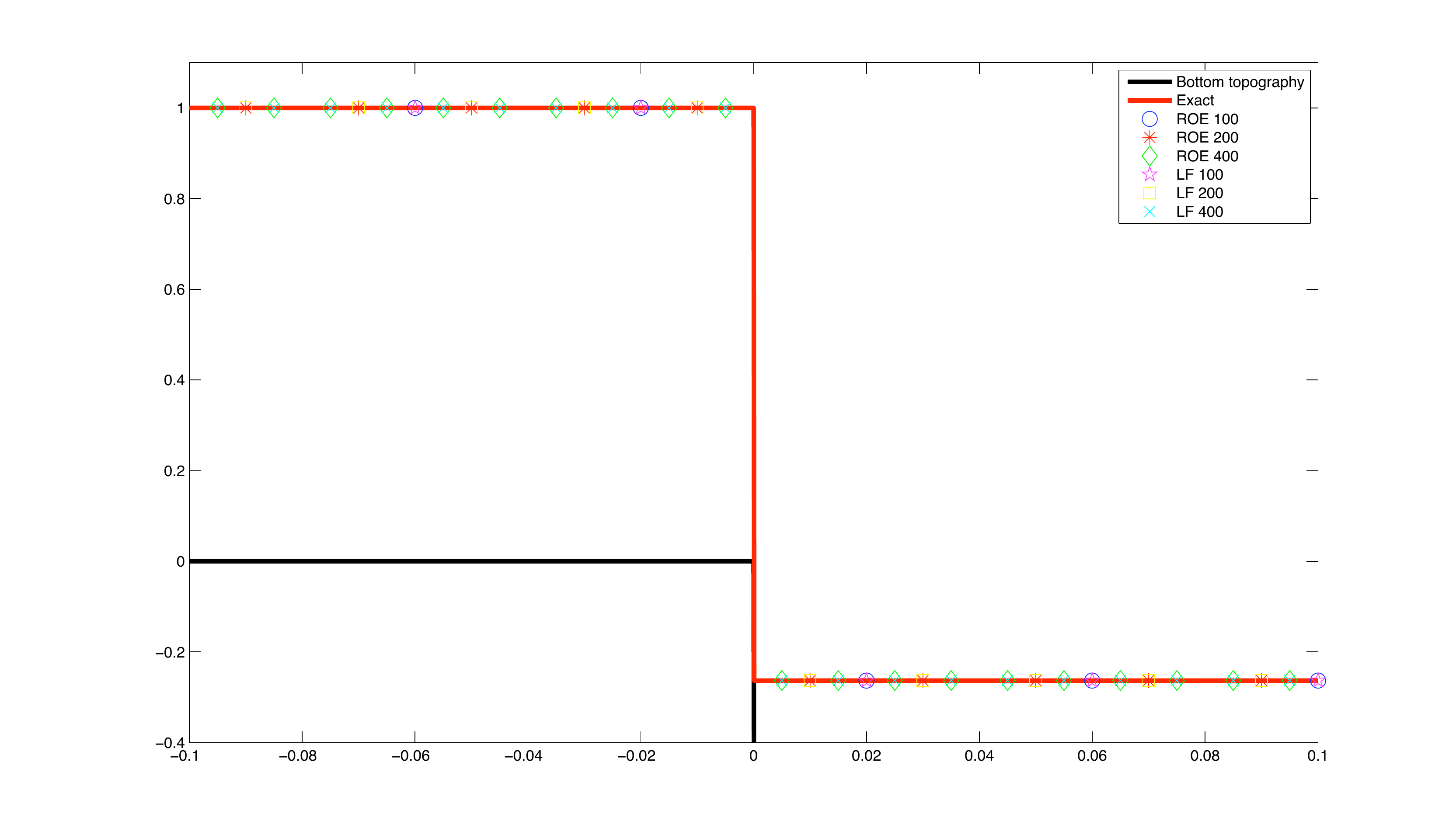}
\end{minipage}}
\subfigure[Bottom topography and free surface (stationary solution): zoom]{\
\begin{minipage}[b]{0.48\textwidth}
\centering
\includegraphics[width=6.5cm,height=6.0cm]{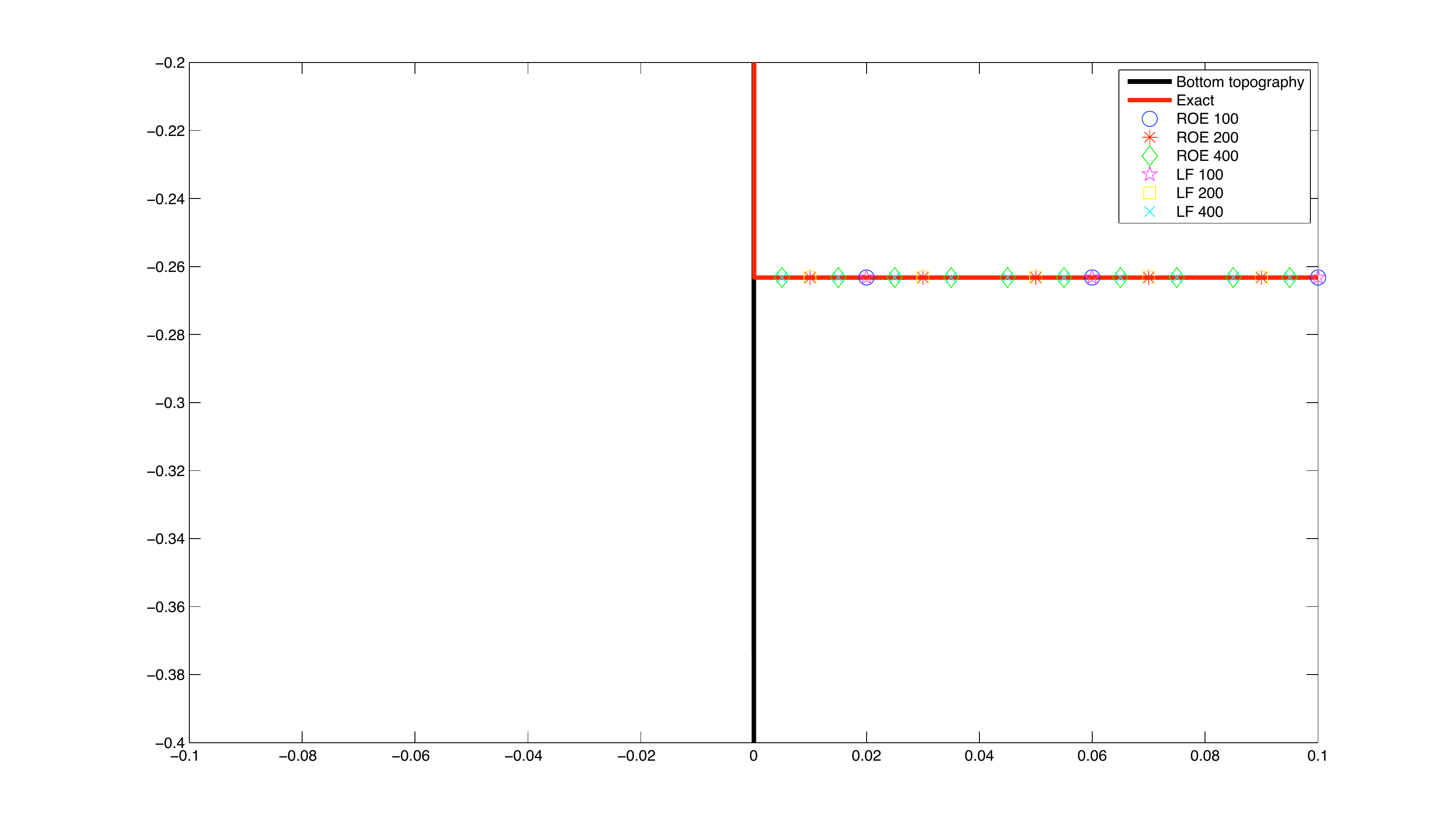}
\end{minipage}}
\caption{Test case \ref{test2.2}. Stationary contact discontinuity: Comparison between the exactly well-balanced Lax-Friedrichs and Roe schemes and the exact solution (free surface).}
\label{test22_fig3}
\end{figure}
\end{center}

Figure \ref{test22_fig3} shows the comparison between  Roe, the modified Lax-Friedrichs scheme, and the exact solution for three meshes with increasing number of cells (100, 200 and 400 cells respectively). As expected, the stationary contact discontinuity is {\sl exactly} captured.

%--------------------------------------------------------------------------------------------------------------

\subsection{Two-layer shallow water system}
\label{2lsw}

\subsubsection{Approximation of internal shocks}
\label{test4.1}

In this section we consider the discretization of the homogeneous two-layer shallow water system (that is, $H = cst.$) by means of Roe and  Lax-Friedrichs schemes. We check that the numerical solutions do not converge to the weak solutions involving shocks even when the same family of paths is used both for the definition of the jump conditions and the construction of the numerical scheme.

We begin with the family of segments for the definition of the jump condition \eqref{RH-2C-1}. In this case:
$$
\aligned
 & \int_0^1 \Phi_{h_1}(s; w_l, w_r)\frac{\partial  \Phi_{h_2} }{\partial s}(s; w_l, w_r)\, ds  = 
 \frac{h_1^l + h_1^r}{2} (h_2^r - h_2^l),\\
 & \int_0^1 \Phi_{h_2}(s; w_l, w_r)
\frac{\partial  \Phi_{h_1}}{ \partial s}(s; w_l, w_r)\,ds
 =   \frac{h_2^l + h_2^r}{2} (h_1^r - h_1^l).
\endaligned
$$

In all of the cases considered here the order of the eigenvalues of the system is:
$$
\lambda_{ext}^- < \lambda_{int}^- <  \lambda_{int}^+ < \lambda_{ext}^+.$$
Moreover
\be
|\lambda_{int}^\pm| << |\lambda_{ext}^\pm| .
\label{eig-ineq}
\ee

We consider first  a Lax-Friedrichs scheme consistent with  the family of segments.  Some easy calculations show that, for the particular choice of the family of segments, the last term
in the modified equation \eqref{modified_2} vanishes.

%We also use the same family of segments to construct the Roe
%linearization (see  \cite{CastroPares_M2AN_roe} for details) used to define  Lax-Friedrichs method \eqref{scheme}-(\ref{DLFM}) :

\begin{center}
\begin{figure}[h!ptb]
\subfigure[Hugoniot curves (projection onto the plane $h_1-q_1$): exact (continuous blue line) and numerical (lines with dots)]{\
\begin{minipage}[b]{0.48\textwidth}
\centering
\includegraphics[width=6.5cm,height=6.0cm]{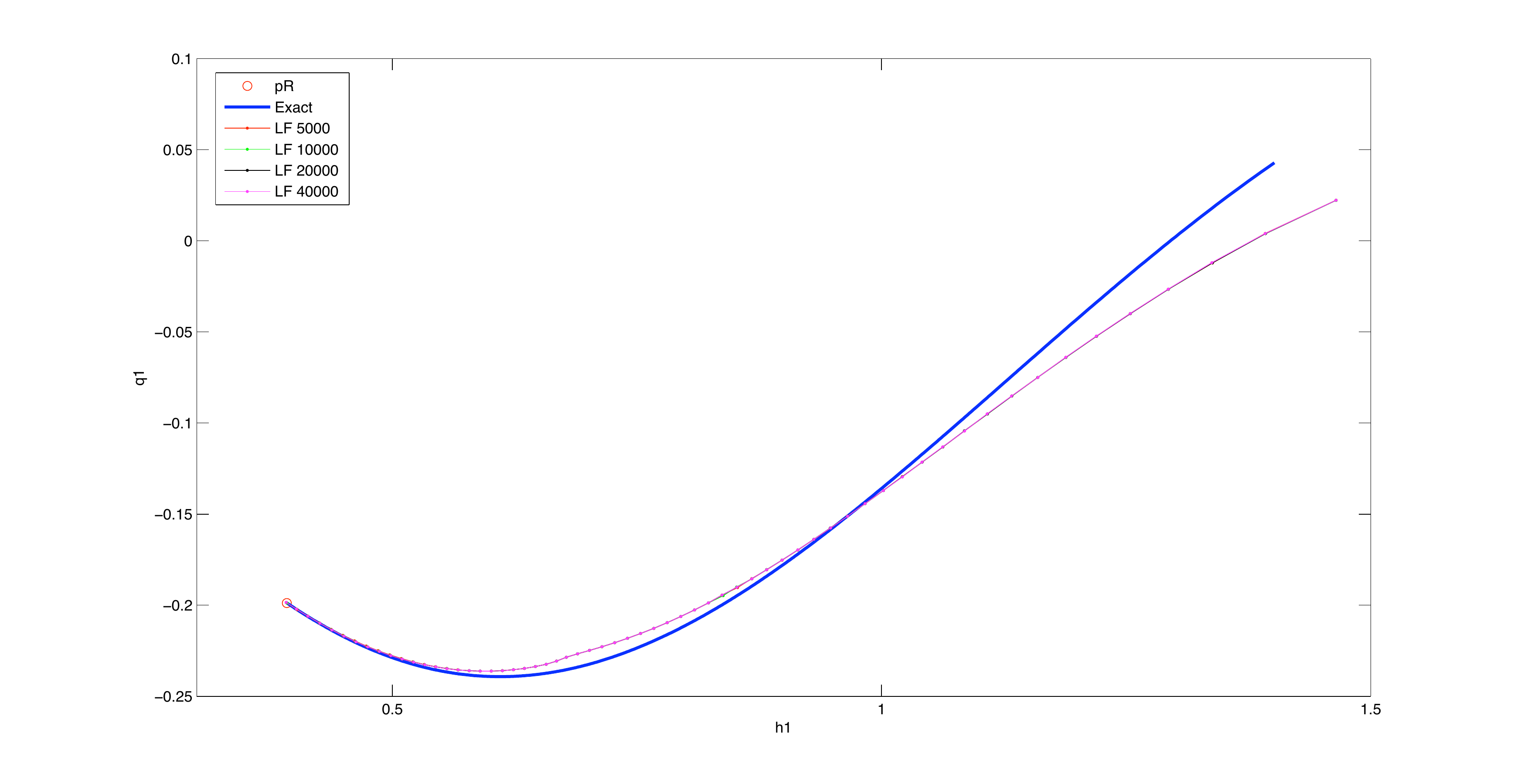}
\end{minipage}}
\subfigure[Hugoniot curves (projection onto the plane $h_2-q_2$): exact (continuous blue line) and numerical (lines with dots)]{\
\begin{minipage}[b]{0.48\textwidth}
\centering
\includegraphics[width=6.5cm,height=6.0cm]{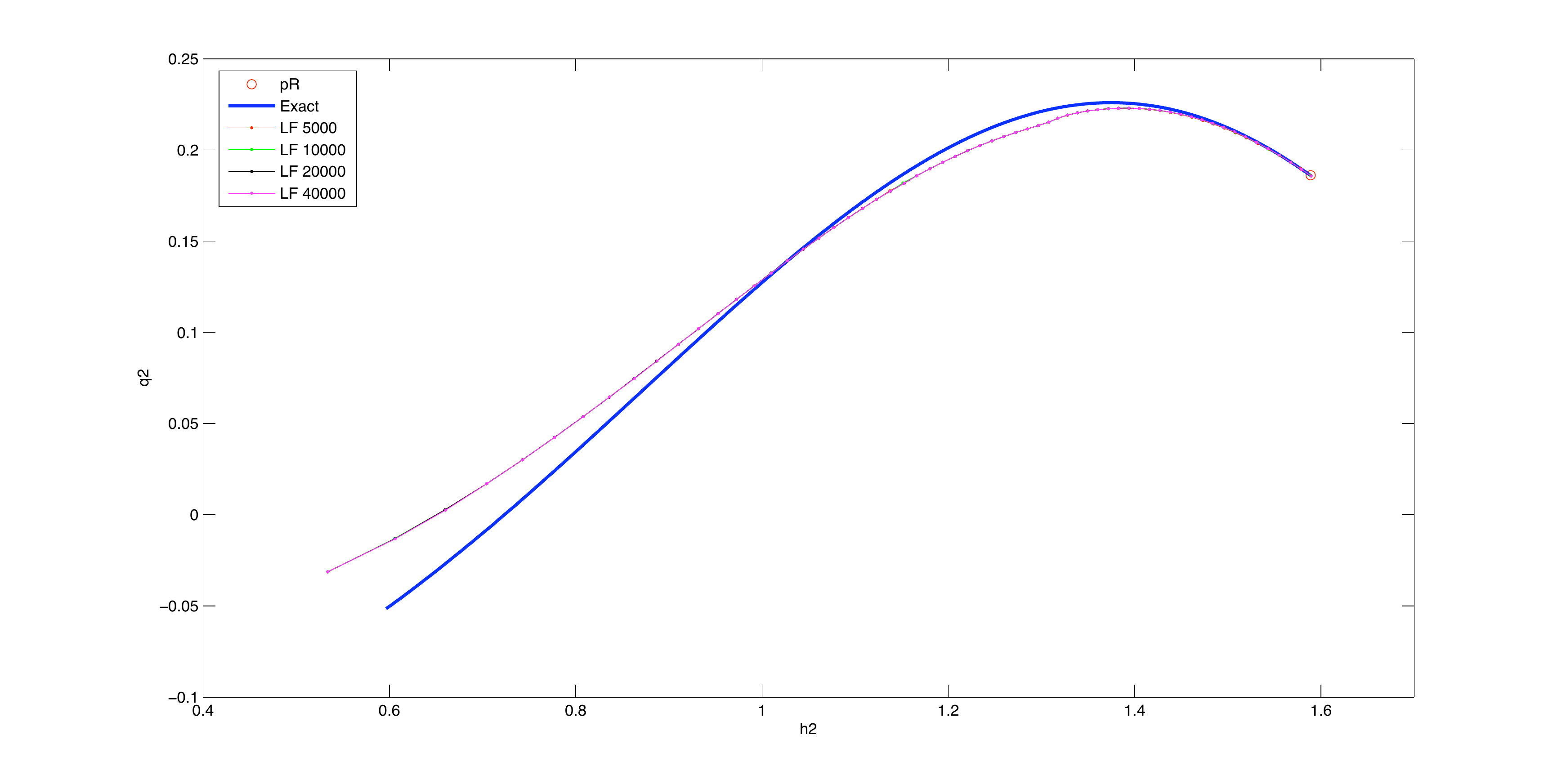}
\end{minipage}}
\caption{Test case \ref{test4.1}. Hugoniot curves: exact (continuous blue line) and numerical (lines with dots).}
\label{test4.1_fig1}
\end{figure}
\end{center}

%\begin{equation}\label{ROEMATSW-2c}
%A^n_{i+1/2} =
%  \left[ \begin{array}{cccc}
%0 & 1 & 0 & 0 \\
%-(u^n_{i+1/2,1})^2 + (c^n_{i+1/2,1})^2 & 2 u^n_{i+1/2,1} &
% (c^n_{i+1/2,1})^2 & 0 \\
%0 & 0 & 0 & 1 \\
%r (c^n_{i+1/2,2})^2 & 0 & -(u^n_{i+1/2,2})^2 + (c^n_{i+1/2,2})^2 & 2 u^n_{i+1/2,2} \\
%\end{array}
%\right].
%\end{equation}
%Here,
%$$
%u^n_{i+1/2,k} = \frac{\sqrt{h^n_{i,k}}u^n_{i,k} +
%                             \sqrt{h^n_{i+1,k}} u^n_{i+1,k}}
%                            {\sqrt{h^n_{i,k}} + \sqrt{h^n_{i+1,k}}},
%\quad  c^n_{i+1/2,k} = \sqrt{g\frac{h^n_{i,k} + h^n_{i+1,k}}{2}},
%k=1,2.
%$$
The goal is to compare the exact  and the numerical Hugoniot curves corresponding to one of the internal characteristic
fields; i.e. the
fields related to the eigenvalues $\lambda_{int}^\pm$).
We proceed as follows:  the state
\begin{equation} \label{WR}
w_r=\left[
\begin{array}{c}
h_1^r \\
q_1^r \\
h_2^r \\
q_2^r
\end{array}
\right]=
\left[
\begin{array}{c}
 0.392034161025472\\
  -0.198826959396196\\
   1.588829011097482  \\
    0.186046955388750
\end{array}
\right]
\end{equation}
is fixed. Then, we compute the Hugoniot curve corresponding to the \lq left \rq states $w_l$ that can be connected with $w_r$ with a 3-shock. To do this, we use the speed of the shock $\xi$ as a parameter and, for each value of $\xi$ we solve the
non-linear system \eqref{RH-2C-1}.  In Figure \ref{test4.1_fig1} we show  the projection (continuous blue line) of the computed Hugoniot curve onto the planes $h_1-q_1$ (left) and $h_2-q_2$ (right), respectively.

Next,  we solve numerically a family of Riemann problems in which the right state is $w_r$,
while $w_l$ runs on the Hugoniot curve.  Using the first divided difference  as a smooth indicator, the speed of propagation 
and the limit states of the shock corresponding to the eigenvalue $\lambda_{int}^+$ is determined in the numerical solution.
These calculations have been performed  by using  four meshes with decreasing step ( $\Delta x= 0.002$, $0.001$, $0.0005$ and $0.00025$).   The numerical Hugoniot curves so obtained are compared with the exact one in Figure  \ref{test4.1_fig1}. Observe that  the numerical solutions converge, but the limit is not a weak solution according to the chosen family of paths. Nevertheless, if $w_r$ and $w_l$ are close enough, both curves are very close. The same behaviour can be observed if the shock speed is close to zero:  this situation corresponds in the figure   to the intersections of  the curves.

\begin{center}
\begin{figure}[h!ptb]
\subfigure[Hugoniot curves (projection onto the plane $h_1-q_1$): exact (continuous red line) and numerical (lines with dots)]{\
\begin{minipage}[b]{0.48\textwidth}
\centering
\includegraphics[width=6.5cm,height=6.0cm]{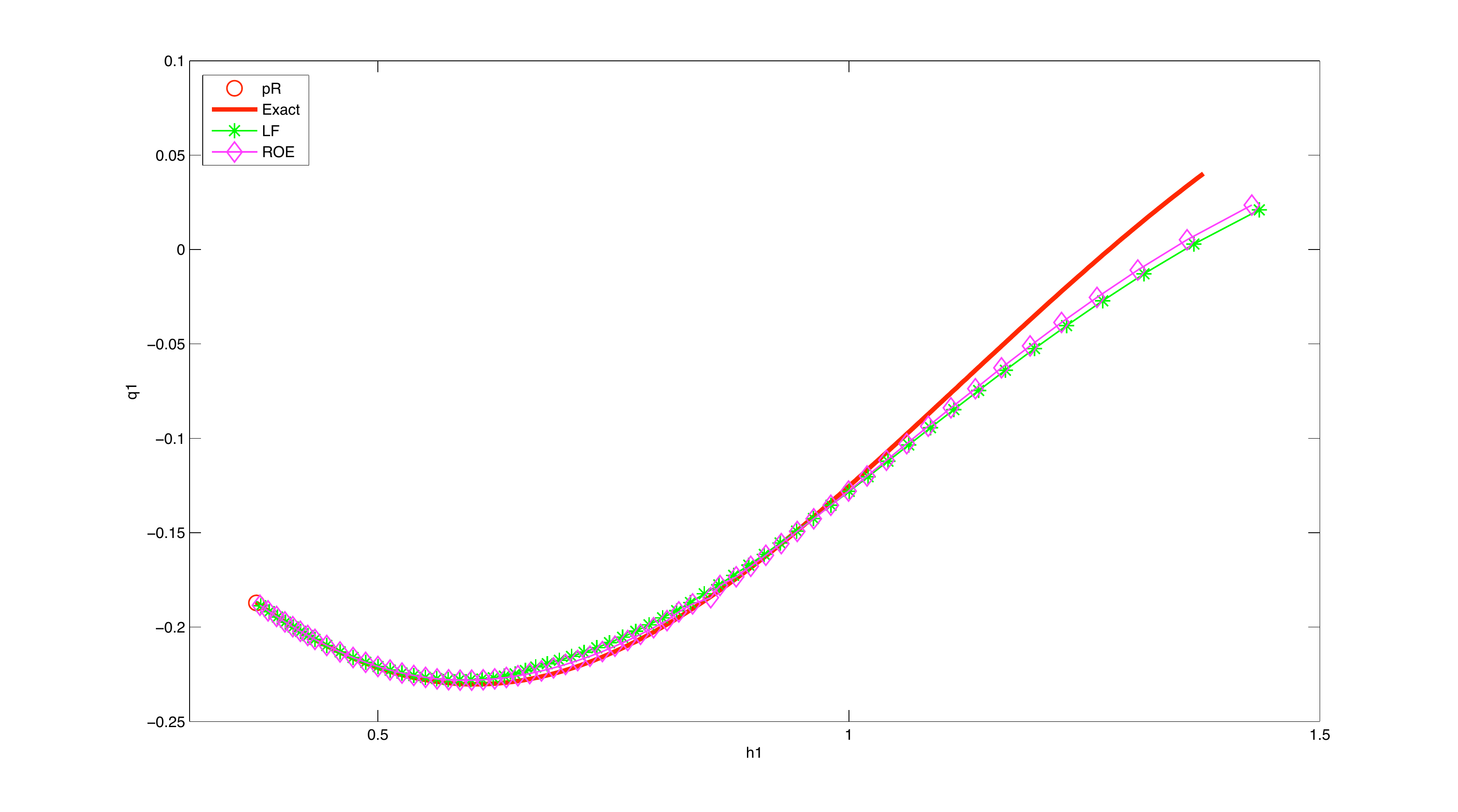}
\end{minipage}}
\subfigure[Hugoniot curves (projection onto the plane $h_2-q_2$): exact (continuous red line) and numerical (lines with dots)]{\
\begin{minipage}[b]{0.48\textwidth}
\centering
\includegraphics[width=6.5cm,height=6.0cm]{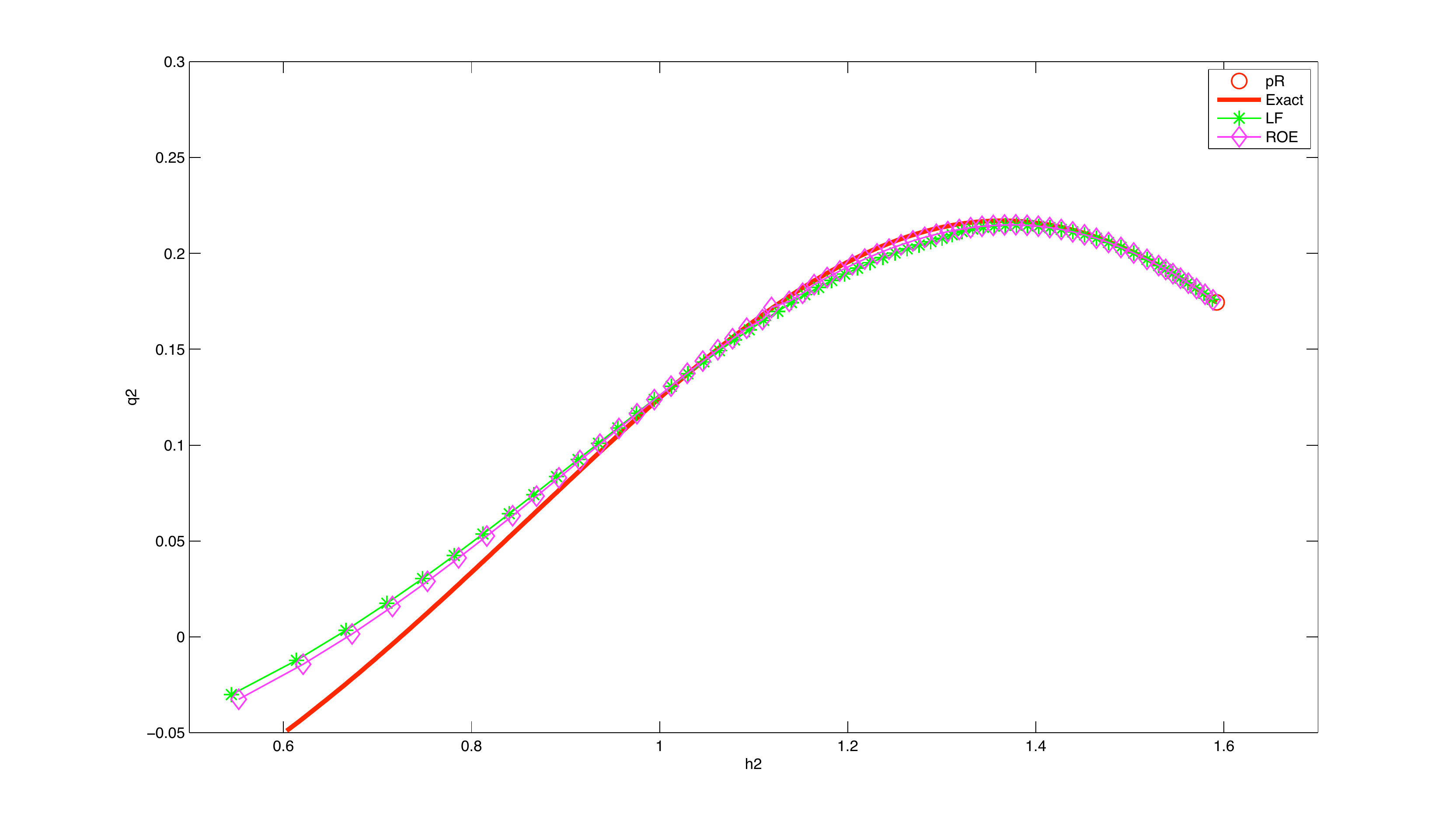}
\end{minipage}}
\caption{Test case \ref{test4.1}. Hugoniot curves: exact (continuous red line) and numerical (lines with dots).}
\label{test4.1_fig2}
\end{figure}
\end{center}

A similar behavior is observed for Roe scheme: the numerical approximations converge but the limit is not a weak solution
according to the family of segments. Nevertheless, as the numerical viscosity of the scheme is smaller than that corresponding
to Lax-Friedrichs, the results are expected to be closer to the exact solution: in Figure  \ref{test4.1_fig2} we compare the
exact Hugoniot curve with those computed with Lax-Friedrichs and Roe schemes using a mesh with 10000 cells.

%-------------------------------------------------------------------------------------------------

\subsubsection{Approximation of external shocks}
\label{test4.2}

Due to the inequality \eqref{eig-ineq}Êthe  CFL condition adjusts the
numerical velocity to the external eigenvalues. As a consequence the effects of the numerical viscosity are much
stronger for internal shocks and thus external shocks are expected to be better captured with
Lax-Friedrichs or Roe schemes.

In order to check this, we proceed as in the previous test case:  the state $w_r$ given by
$$
w_r=\left[
\begin{array}{c}
0.257381469591567 \\
0.444901654188681 \\
0.110306344093418 \\
0.190672137450279
\end{array}
\right]
$$
is fixed, and then the Hugoniot curves corresponding to the \lq left \rq states ($w_l$) that can be connected with $w_r$ with a shock related to the eigenvalue $\lambda^-_{ext}$  are computed by solving the non-linear system \eqref{RH-2C-1} using $\xi$ as a parameter.

\begin{center}
\begin{figure}[h!ptb]
\subfigure[Hugoniot curves (projection onto the  plane $h_1-q_1$): exact (continuous red line) and numerical (lines with dots)]{\
\begin{minipage}[b]{0.48\textwidth}
\centering
\includegraphics[width=6.5cm,height=6.0cm]{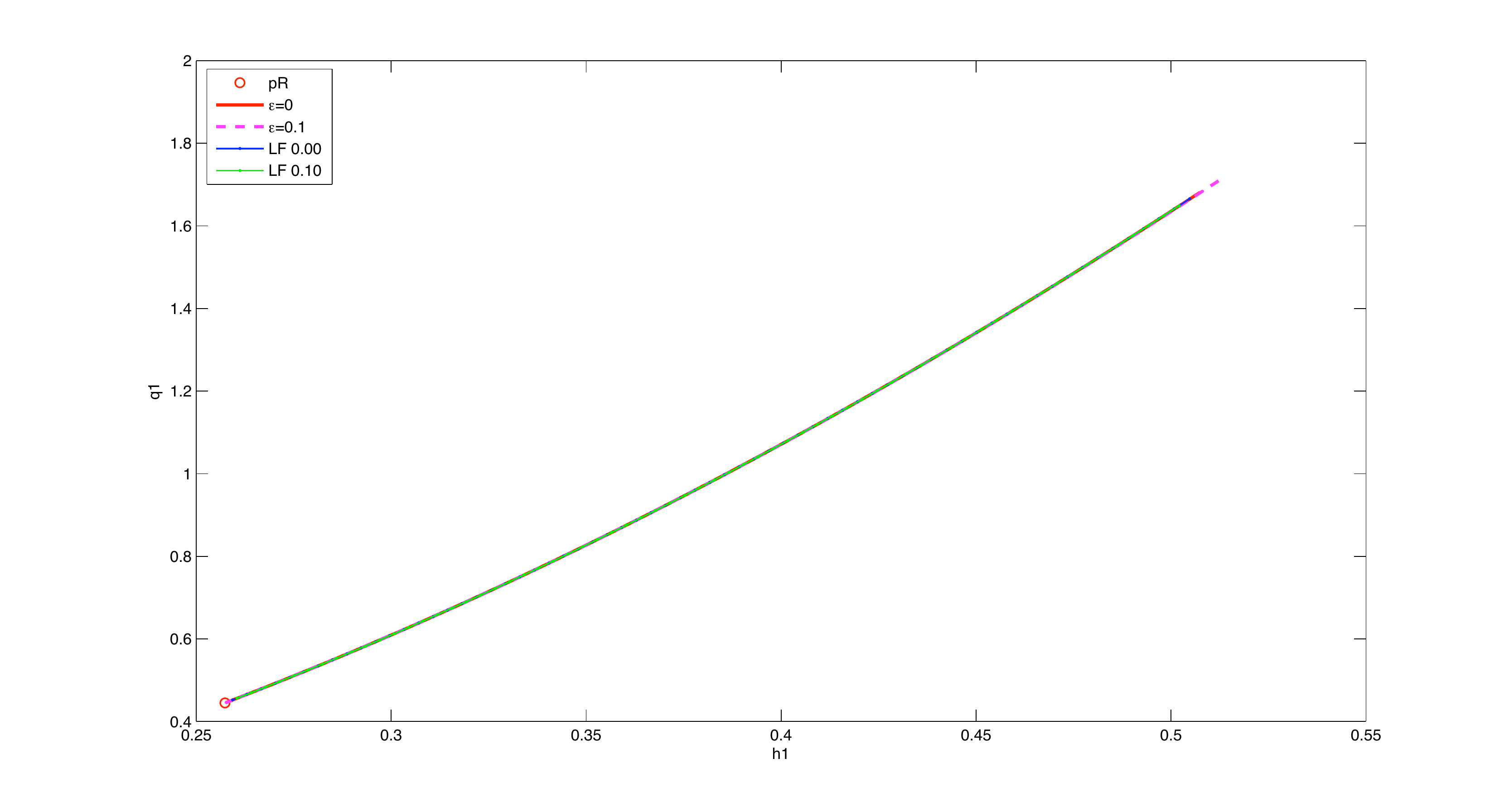}
\end{minipage}}
\subfigure[Hugoniot curves (zoom of the projection onto the plane $h_1-q_1$): exact (continuous red line) and numerical (lines with dots)]{\
\begin{minipage}[b]{0.48\textwidth}
\centering
\includegraphics[width=6.5cm,height=6.0cm]{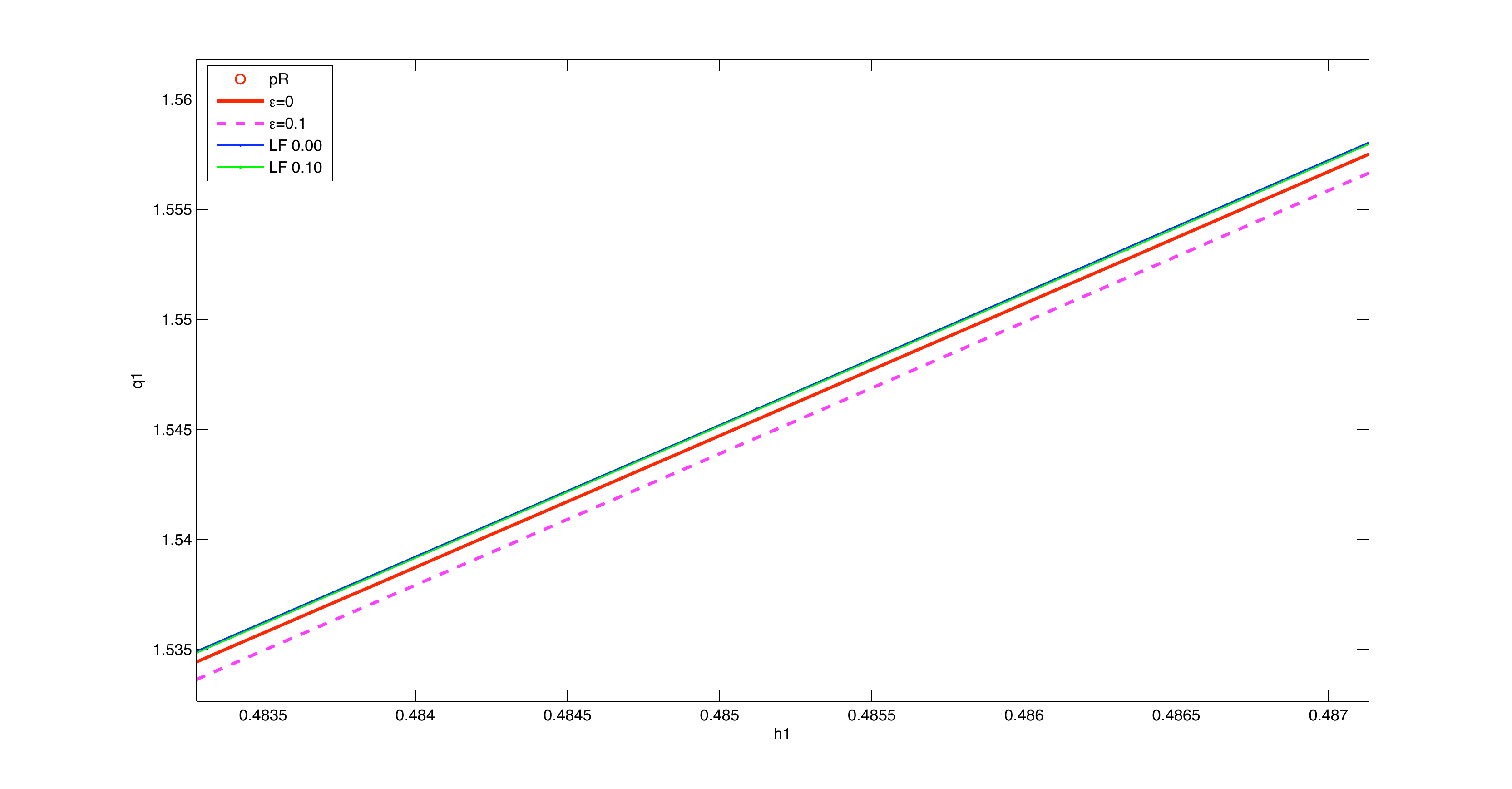}
\end{minipage}}
\subfigure[Hugoniot curves (projection onto the  plane $h_2-q_2$): exact (continuous red line) and numerical (lines with dots)]{\
\begin{minipage}[b]{0.48\textwidth}
\centering
\includegraphics[width=6.5cm,height=6.0cm]{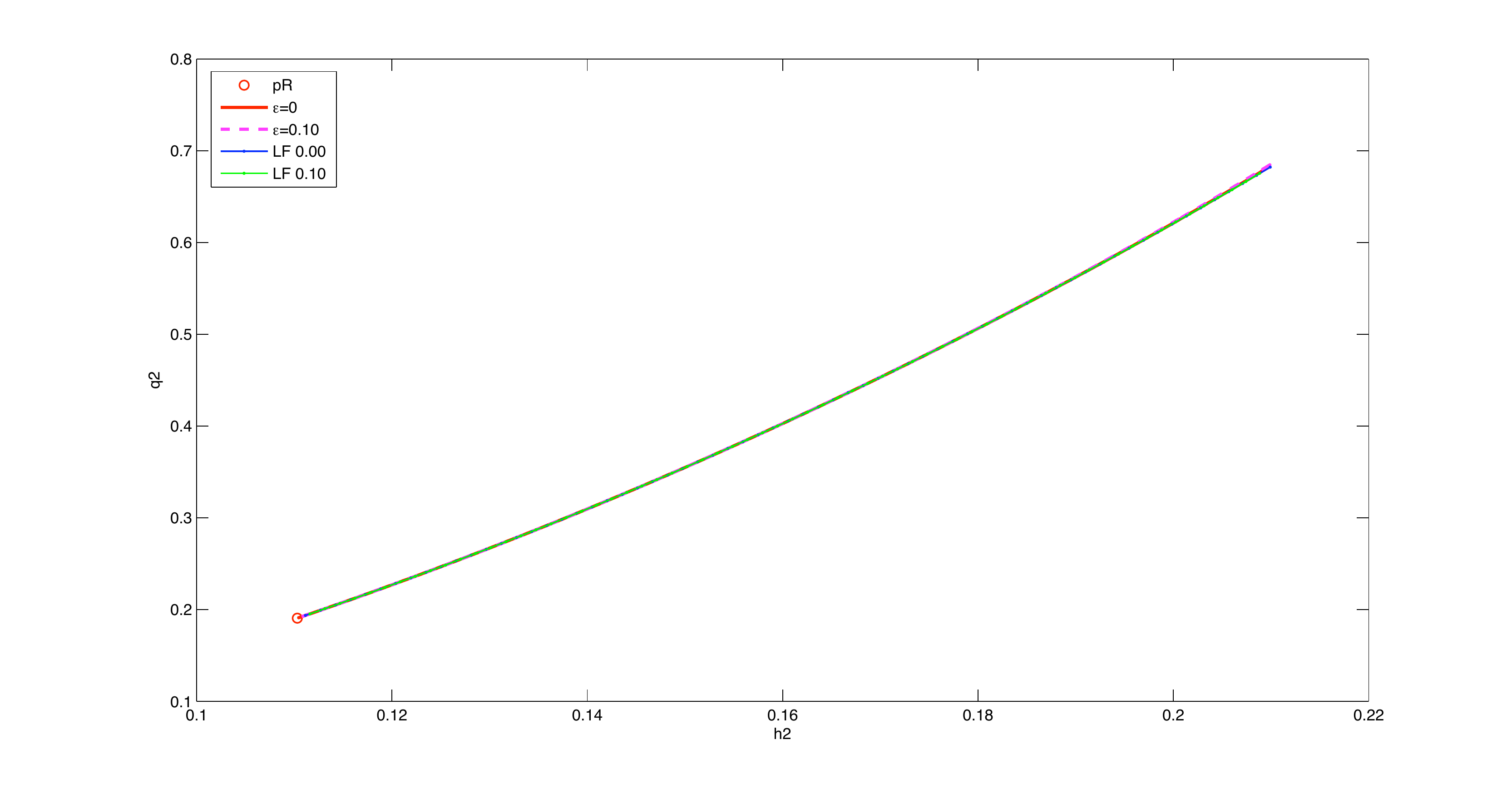}
\end{minipage}}
\subfigure[Hugoniot curves (zoom of the projection onto the plane $h_2-q_2$): exact (continuous red line) and numerical (lines with dots)]{\
\begin{minipage}[b]{0.48\textwidth}
\centering
\includegraphics[width=6.5cm,height=6.0cm]{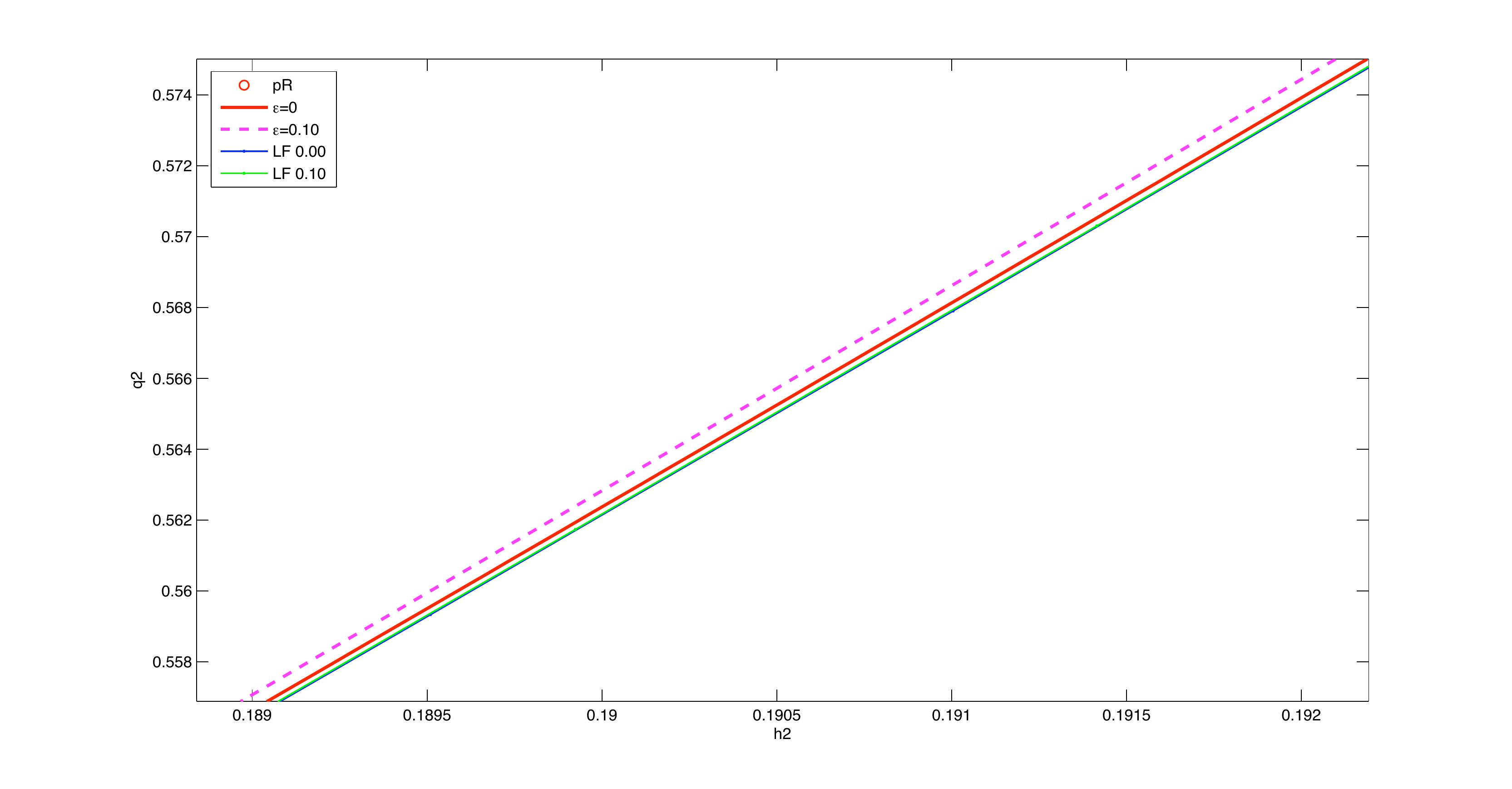}
\end{minipage}}
\caption{Test case \ref{test4.2}. Hugoniot curves for Lax-Friedrichs and Roe  schemes: exact (continuous red line) and numerical (lines with dots).}
\label{test4.2_fig1}
\end{figure}
\end{center}

In Figure \ref{test4.2_fig1} we compare the
exact Hugoniot curve with those computed with Lax-Friedrichs and Roe scheme using a mesh with 2000 cells.
Note how the curves are  now much closer to each other than they were in the previous test case.

%-------------------------------------------------------------------------------------------------------------

\subsection{Influence of the family of paths}
\label{test4.3}

In this test we study the influence of small changes in the family of paths both in the weak  and the numerical solutions.
We consider now the jump conditions related to a different family of paths  $\Phi^\epsilon_{h_i}(s; w_l, w_r)$, $i=1,2$.
These curves are chosen so that:
$$
\aligned
& h_1  =  \Phi^\epsilon_{h_1}(s; w_l, w_r), \quad s \in [0,1];\\
& h_2  =  \Phi^\epsilon_{h_2}(s; w_l, w_r), \quad s \in [0,1];\\
\endaligned
$$
is a parameterization of the curve
$$
h_2= h_2^l+ \left( \frac{h_1-h_1^l}{h_1^r-h_1^l}+\epsilon \frac{(h_1^r)^2-(h_1^l)^2}{(h_1^r)^2-(h_1^l)^2}\right) \frac{h_2^r-h_2^l}{1+\epsilon}.
$$

The jump conditions are now  \eqref{RH-2C-1}, with
$$
 \int_0^1 \Big( \Phi_{h_1} \, \frac{\partial \Phi_{h_2} }{\partial s}\Big) (s; w_l, w_r)\, ds  = 
    \frac{(h_1^l)^2(3+4\epsilon) +2(3+2\epsilon)h_1^l h_1^r +(3+4\epsilon)(h_1^r)^2 }{6(1+\epsilon)(h_1^r+h_1^l)},
$$
$$
 \int_0^1 \Big( \Phi_{h_2}
\frac{\partial  \Phi_{h_1}}{ \partial s}\Big) (s; w_l, w_r)\,ds  = 
 \frac{h_1^r\left( (3+4\epsilon)h_2^l+(3+2\epsilon)h_2^r\right) + h_1^l\left( (3+2\epsilon) h_2^l+(3+4\epsilon)h_2^l\right)}{6(1+\epsilon)(h_1^r+h_1^l)}.
$$
Observe  that the jump conditions of the previous tests are recovered for  $\epsilon=0$.

\begin{center}
\begin{figure}[h!ptb]
\subfigure[Hugoniot curves (projection onto the plane $h_1-q_1$)]{\
\begin{minipage}[b]{0.48\textwidth}
\centering
\includegraphics[width=6.5cm,height=6.0cm]{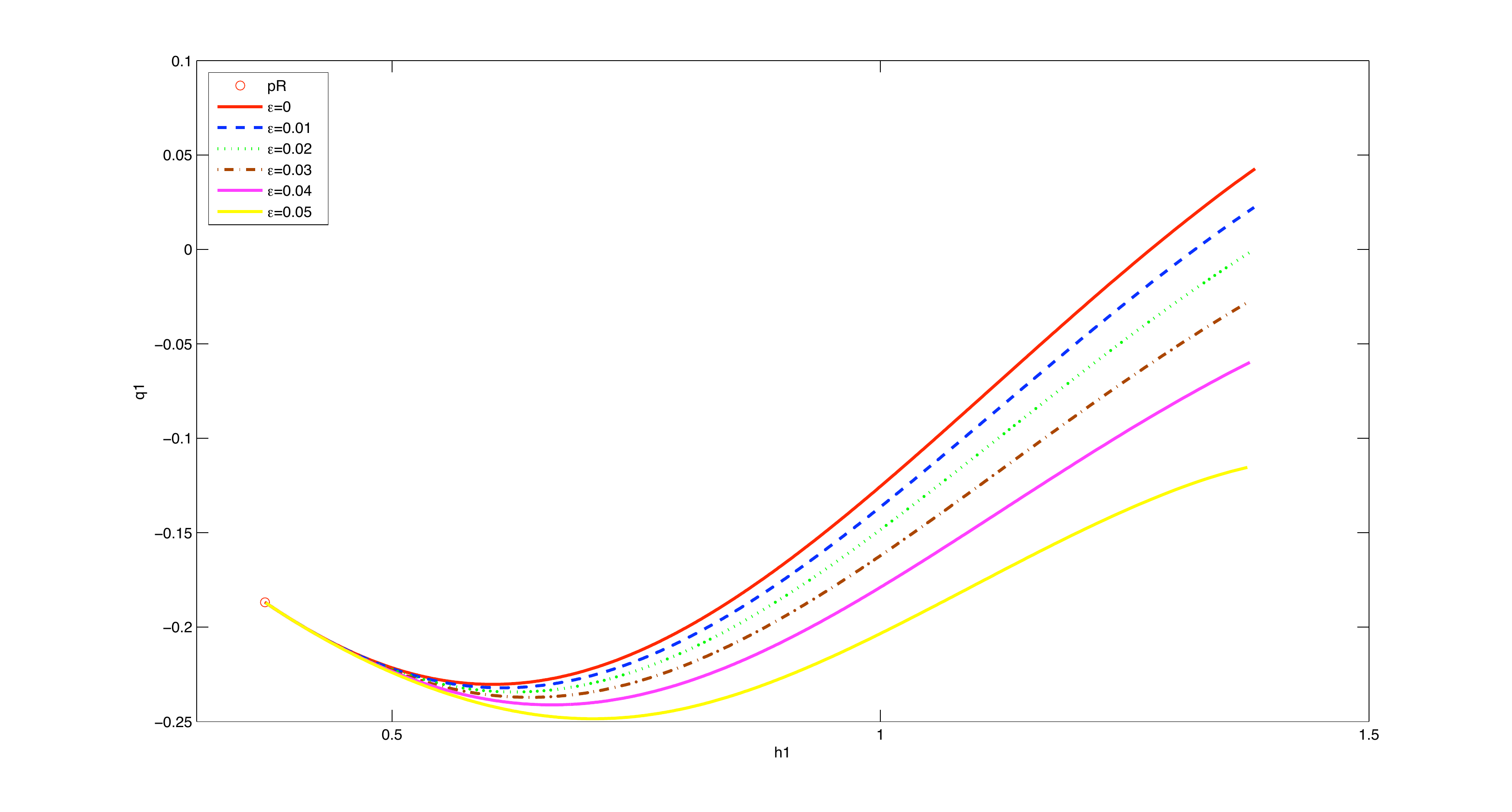}
\end{minipage}}
\subfigure[Hugoniot curves (projection onto the plane $h_2-q_2$)]{\
\begin{minipage}[b]{0.48\textwidth}
\centering
\includegraphics[width=6.5cm,height=6.0cm]{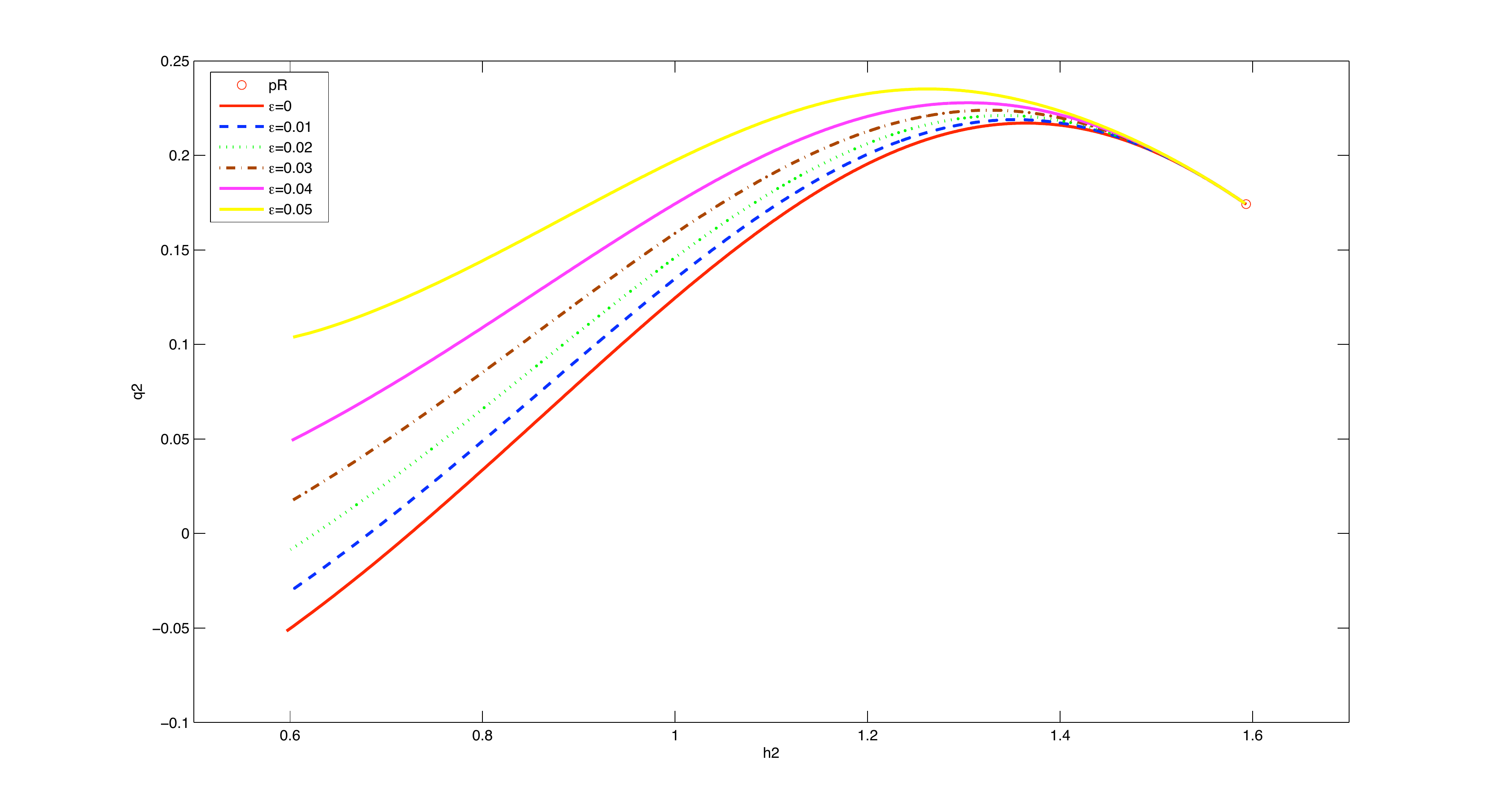}
\end{minipage}}
\caption{Test case \ref{test4.3}. Hugoniot curves: \mbox{$\epsilon \in \left\{0.0,\ 0.01,\ 0.02, \ 0.03,\ 0.04, \ 0.05 \right\}$}.}
\label{test4.3_fig1}
\end{figure}
\end{center}

As in Section~\ref{test4.1},  the state $w_r$ given by \eqref{WR} is fixed and the exact Hugoniot curves are computed for different values of $\epsilon$. In Figure \ref{test4.3_fig1}  the projections of the exact Hugoniot curves onto the planes $h_1-q_1$ (left) and $h_2-q_2$(right) for the values $\epsilon$ = 0.00, 0.01, 0.02,  0.03, 0.04, 0.05, are shown.

Next, we consider the  Lax-Friedrichs schemes related to the new choice of the family of paths. An easy calculation shows again  that, for this new family of paths, the last term of the modified equation  \eqref{modified_2} also vanishes. As a consequence, the second order modified equation is independent of $\epsilon$: it coincides with the one corresponding to the family of
segments.

\begin{center}
\begin{figure}[h!ptb]
\subfigure[Hugoniot curves (projection onto the  plane $h_1-q_1$)]{\
\begin{minipage}[b]{0.48\textwidth}
\centering
\includegraphics[width=6.5cm,height=6.0cm]{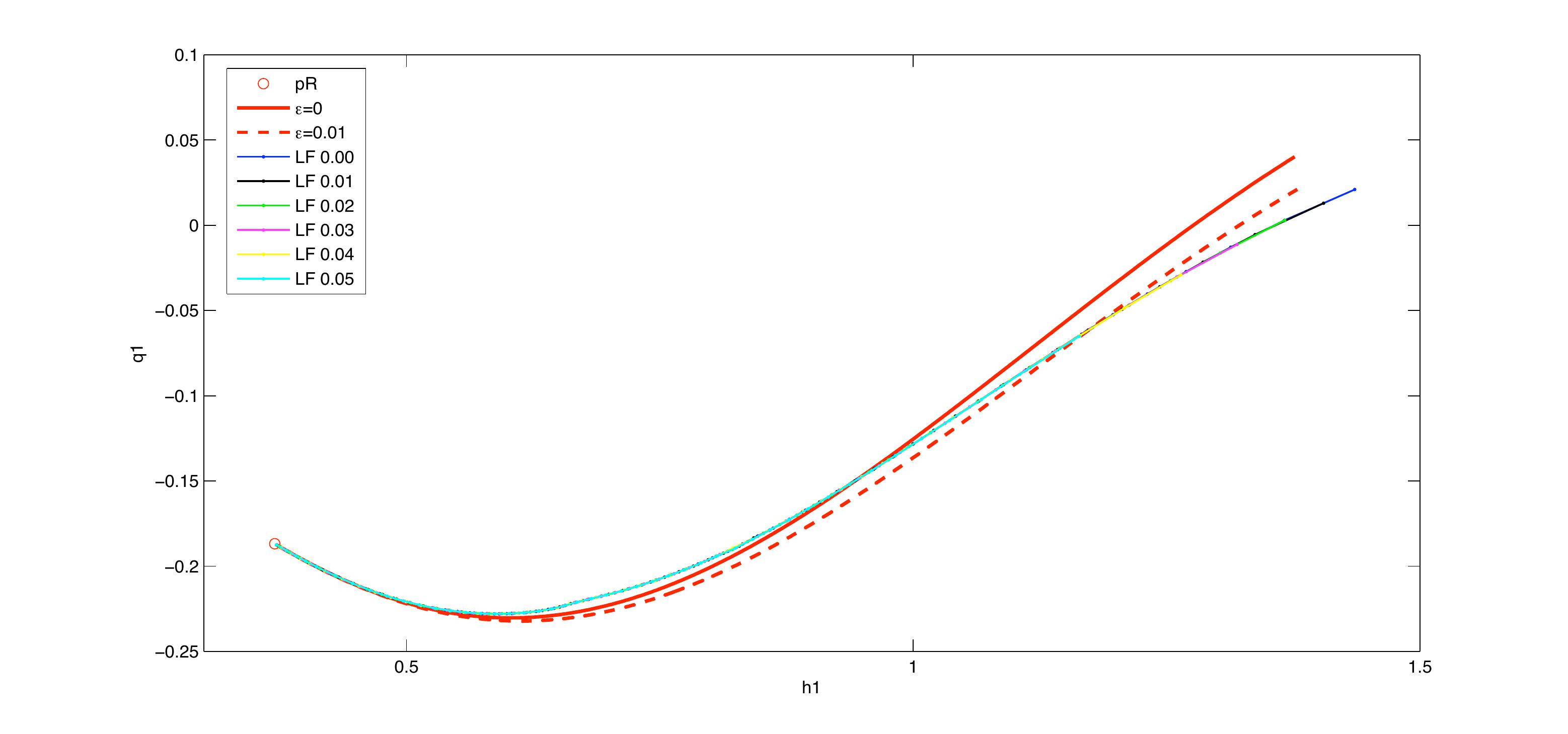}
\end{minipage}}
\subfigure[Hugoniot curves (projection onto the plane $h_2-q_2$)]{\
\begin{minipage}[b]{0.48\textwidth}
\centering
\includegraphics[width=6.5cm,height=6.0cm]{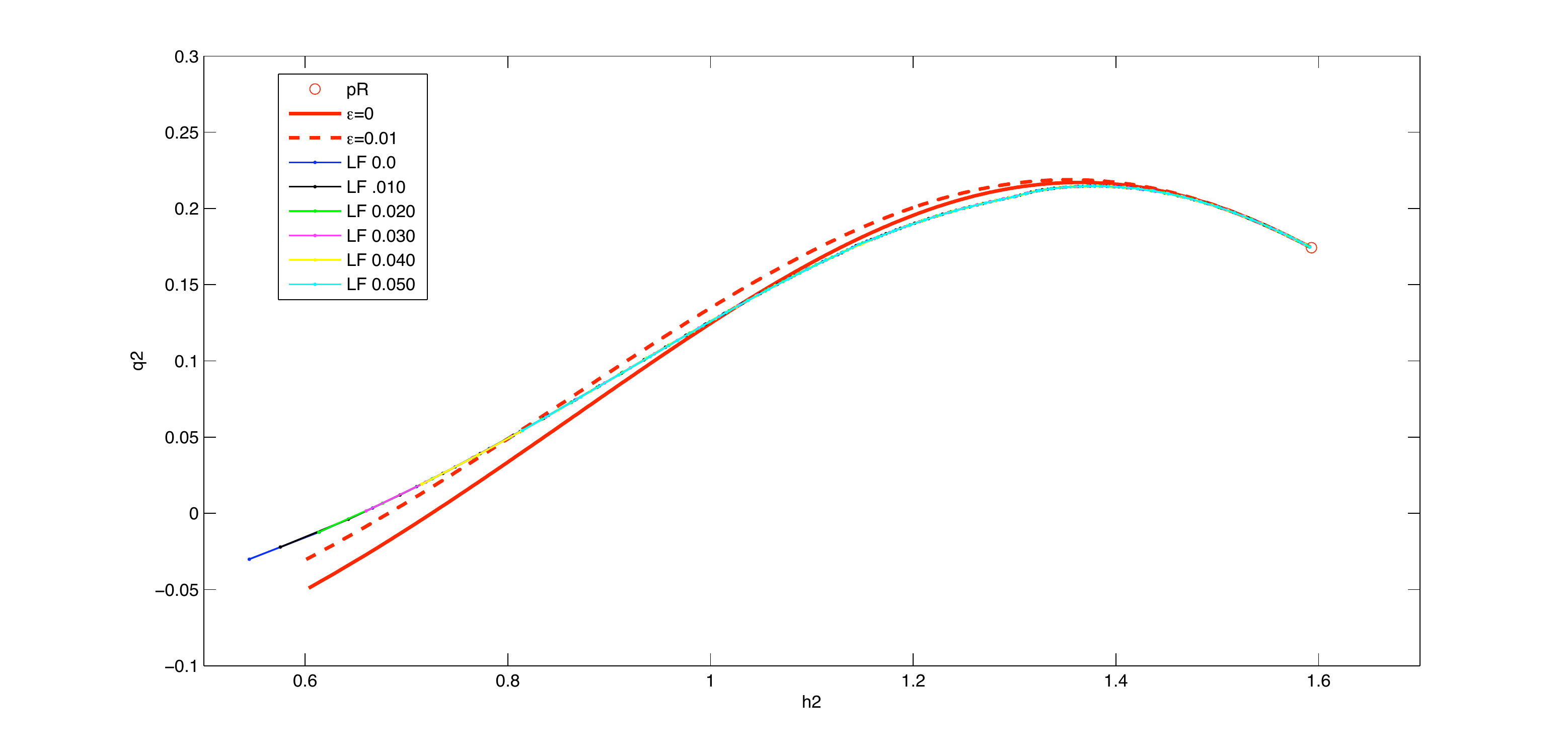}
\end{minipage}}
\caption{Test case \ref{test4.3}. Hugoniot curves for Lax-Friedrichs scheme with $\epsilon \in \left\{0.0,\ 0.01,\ 0.02, \ 0.03,\ 0.04, \ 0.05 \right\}$: exact (red lines) and numerical (lines with dots). }
\label{test4.3_fig2}
\end{figure}
\end{center}

We proceed as in the previous test case to compute the Hugoniot curves corresponding to the numerical scheme.
In Figure \ref{test4.3_fig2} the curves obtained using a mesh with step $\Delta x =0.001$ for $\epsilon$ =   0.00, 0.01, 0.02,
0.03, 0.04, 0.05,  are depicted. Observe that all of the curves coincide:  the Hugoniot curves  obtained with the Lax-Friedrichs scheme for the different values of   $\epsilon$  are reparameterizations of the same curve. This fact agrees with the fact
that the second order modified equation is independent of $\epsilon$.

\begin{center}
\begin{figure}[h!ptb]
\subfigure[3-shock hugoniot curves (projection onto the plane $h_1-q_1$)]{\
\begin{minipage}[b]{0.48\textwidth}
\centering
\includegraphics[width=6.5cm,height=6.0cm]{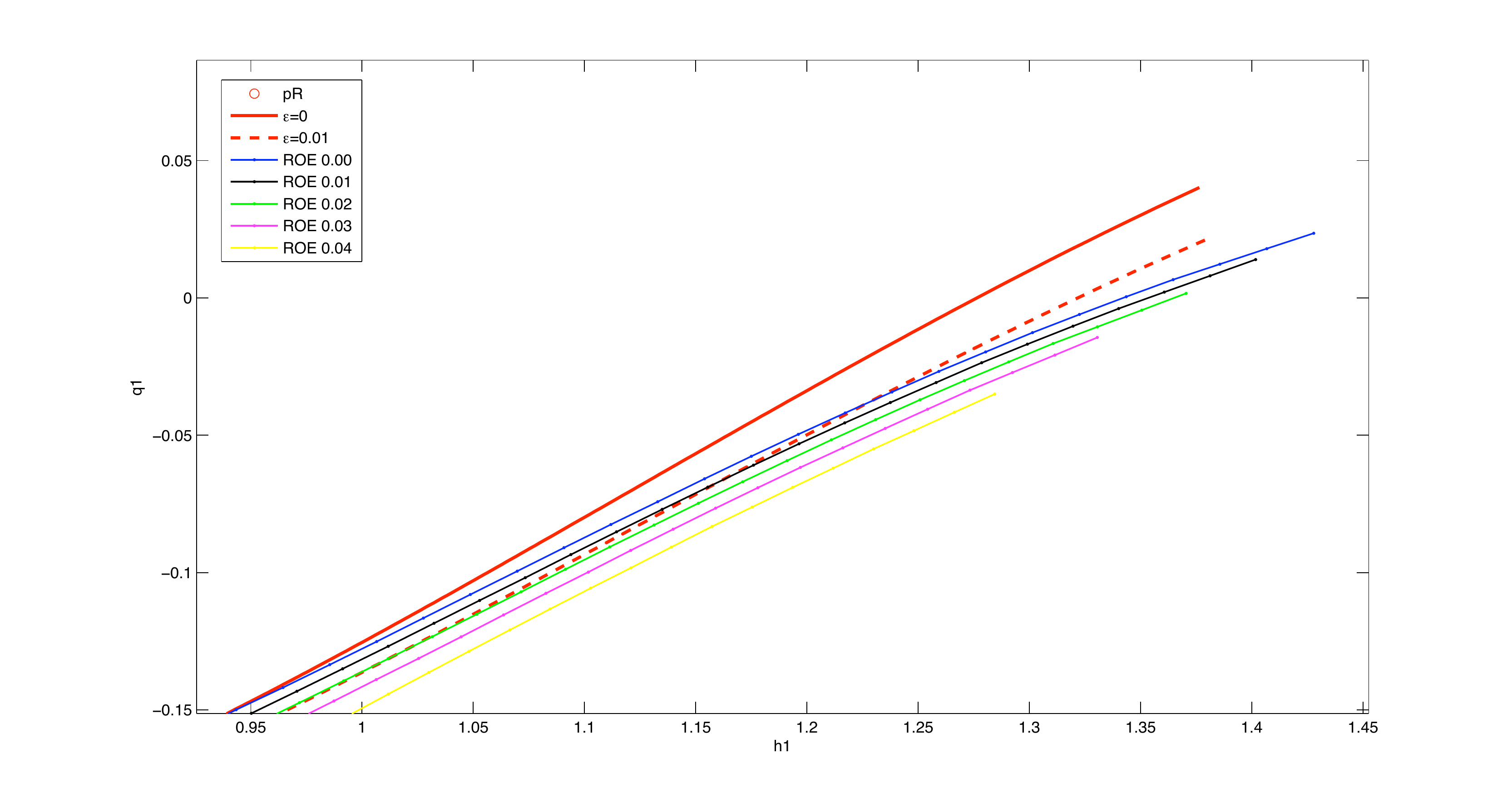}
\end{minipage}}
\subfigure[3-shok hugoniot curves (projection onto the  plane $h_2-q_2$)]{\
\begin{minipage}[b]{0.48\textwidth}
\centering
\includegraphics[width=6.5cm,height=6.0cm]{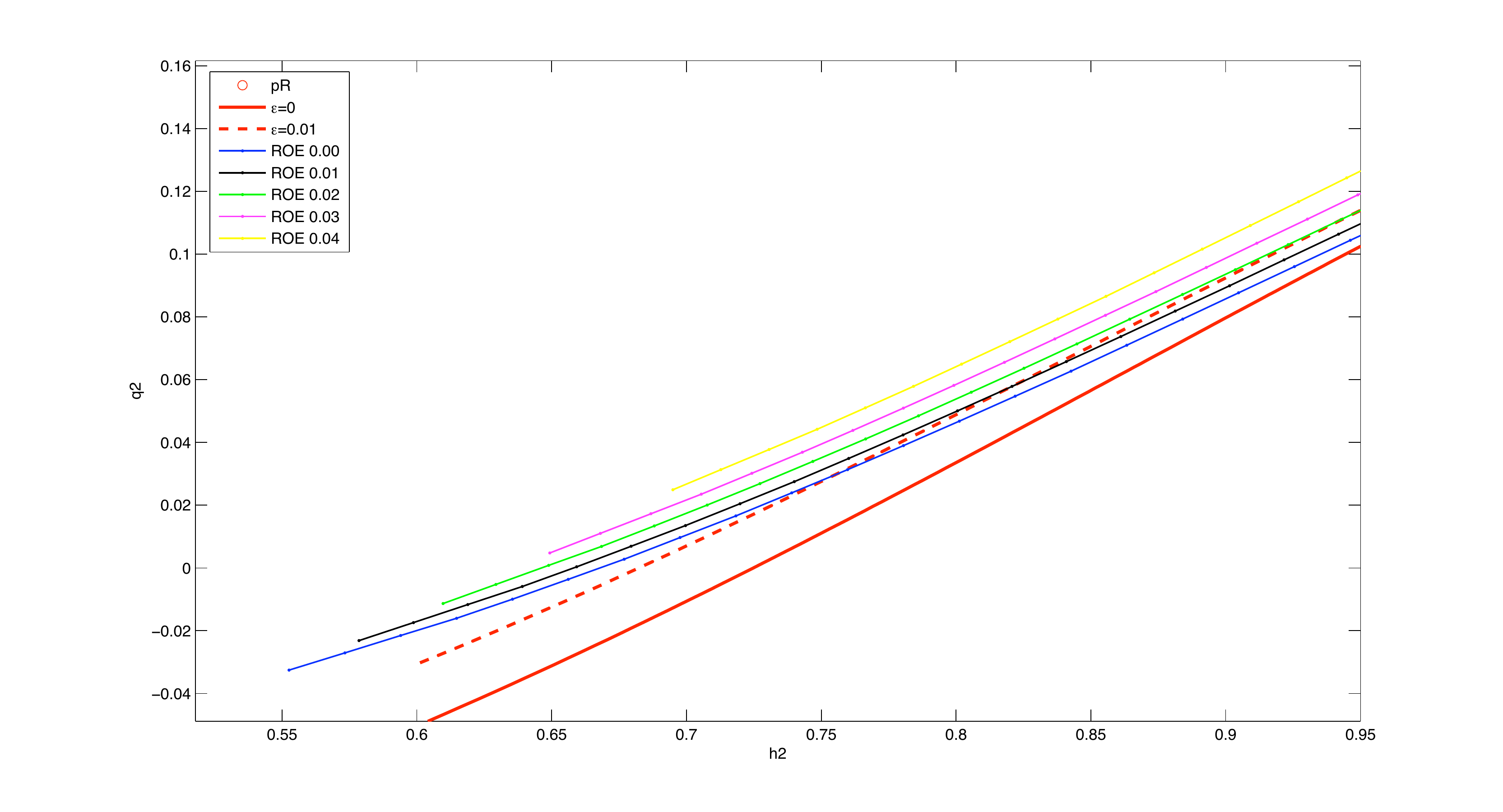}
\end{minipage}}
\caption{Test case \ref{test4.3}. 3-shock hugoniot curves (zoom) for Roe scheme with $\epsilon \in \left\{0.0,\ 0.01,\ 0.02, \ 0.03,\ 0.04, \ 0.05 \right\}$: exact (red lines) and numerical (lines with dots). }
\label{test4.3_fig3}
\end{figure}
\end{center}

Figure \ref{test4.3_fig3} shows the projections of the Hugoniot curve corresponding to  Roe scheme for  $\epsilon$ =  0.00, 0.01, 0.02, 0.03, 0.04, 0.05 (line with dots) onto the planes $h_1-q_1$ and $h_2-q_2$. Note that, in this case, the obtained curves depend on  $\epsilon$, but again the numerical solutions do not converge to the weak solutions
corresponding to the chosen family of paths.  Nevertheless, as in the previous test case,  if $w_r$ and $w_l$ are close enough or the shock speed is close to zero, the exact and the numerical Hugoniot curves are also close.

%*****************************************************************************************
\newpage 

\section{Conclusions}

When a hyperbolic system with nonconservative products and genuinely nonlinear fields
is discretized, in order to be sure that the numerical approximations
converge to a function which is a classical solution where it is
smooth and whose discontinuities are in  good  agreement with the
physics of the problem, the following steps should be taken \cite{LeFloch1,HouLeFloch}: 
{\begin{itemize}
\item First, choose a regularization of the system which is consistent with the physics of the problem.

\item Next, determine the DLM family of paths consistent with this regularization.

\item Finally, design a numerical scheme whose solutions converge to 
weak solutions associated with this family of paths.
\end{itemize} 
In practice, this strategy may be difficult 
to follow, since the actual calculation of a family of paths requires
calculating regularized shock profiles associated with the given regularization. 
On the other hand, the convergence of the numerical solutions to the correct
weak solutions is known for the Glimm scheme and the front tracking method \cite{PLF-Liu}, 
only; 
the implementation of these methods can be time consuming since they require the explicit 
knowledge of the corresponding Riemann solver. 
}

In fact, when the nonconservative model under consideration is a simplified version of 
a more complex (but conservative) model ---as is the case for the two-layer shallow water system--- 
the above strategy may end up being more costly than than solving
directly the more complex one.
In these cases, the use of a numerical strategy based on a direct discretization
of the nonconservative system by means of finite difference scheme which is formally 
path-consistent is advisable and 
may have the following advantages:
\begin{itemize}
\item The numerical solutions is formally consistent
with the definition of the nonconservative product in the sense of Dal~Maso, LeFloch, and Murat and, in turn, 
in the special case that the system admits a conservative subsystem, the numerical scheme 
is conservative for that subsystem in the sense of Lax.

\item The approximations of the shocks provided by the schemes are consistent
with a regularization of the system with higher-order terms that
vanish as $\Delta x$ tends to 0. 
(Obviously, the main drawback is that this  regularization depends
on the chosen family of paths and on the numerical scheme itself. This is issue dealt with in the present 
paper.) 

\item  As originally pointed out by Hou and LeFloch \cite{HouLeFloch} in the (simpler) 
case of scalar hyperbolic equations, 
the convergence error, measured in terms of our convergence error measure  
or in terms of the Hugoniot curves, is noticeable for very fine meshes,
for discontinuities of great amplitude, and/or for large-time simulations, only.

\item This strategy is extendable to high-order methods or to multidimensional problems, as 
developed, together with collaborators,
by Coquel \cite{BC,BerthonCoquel,CC} and Pares \cite{CGP,CGP07}.

\end{itemize}

The convergence error should also be compared with the experimental error. 
In the case of the two-layer shallow water system, the shocks
captured by Roe scheme and the family of straightlines have been found \cite{Cetal}
to be in good agreement with the experimental measurements of internal bores in the Strait of
Gibraltar, despite of the simplicity of the family of paths.

In certain special situations, the convergence error measure is found to vanish identically.
This is the case of systems whose nonconservative product is associated with 
a linearly degenerate field: for schemes that are formally consistent with
a family of paths satisfying the condition (R1), 
then all of the discontinuities are correctly approximated and the scheme does converge to exact solutions. 
The discussion of linearly degenerate fields associated with nonconservative products was discussed 
earlier in \cite{LeFloch1,LeFlochThanh} from the theoretical standpoint and 
in \cite{Bouchut} from the numerical standpoint. 
This problem may also exhibit an additional difficulty, the resonance problem,
if one of the eigenvalues of the Jacobian matrix vanishes, and
weak solutions may not be uniquely determined by their initial data, so that 
the limiting numerical solutions may depend both on the family
of paths and the numerical scheme itself.

%*****************************************************************************************

\section*{Acknowledgments}

This research was partially supported by the Spanish Government Research project MTM2006-08075, 
by the A.N.R. (Agence Nationale de la Recherche) through the grant 06-2-134423, 
and by the Centre National de la Recherche Scientifique (CNRS). The authors are grateful 
to Fr\'ed\'eric Coquel for discussions on the matter of this paper. 

%*****************************************************************************************

\newcommand \auth {\textsc}

\end{document}